\documentclass[12pt,reqno]{amsart}

\RequirePackage{pict2e}
\usepackage{amssymb,amscd,url,amsmath,amscd}
\usepackage{xcolor}    
\usepackage{bbm}       
\usepackage{extarrows}
\usepackage{constants} 
\usepackage{stackrel}
\usepackage{tikz}
\usetikzlibrary{arrows.meta}
\usepackage{centernot}
\usepackage{upref}     
\usepackage{mathtools}

\RequirePackage[pdfa, psdextra]{hyperref}  
\hypersetup{
    colorlinks=true,
    linkcolor=blue,
    urlcolor=cyan,
    }
\providecommand{\texorpdfstring}[2]{#1}  

\begin{document}

\allowdisplaybreaks


\title[Propagation of Zariski Dense Orbits]
{Propagation of Zariski Dense Orbits}
\date{\today}

\author[H. Pasten]{Hector Pasten}
\email{hpasten@gmail.com}
\address{Facultad de Mathm{\'a}ticas, Pontificia Universidad Católica de Chile, Vicu{\~n}a Mackenna 4860, Macul, Chile}

\author[J.H. Silverman]{Joseph H. Silverman}
\email{joseph\_silverman@brown.edu}
\address{Department of Mathematics, Box 1917
  Brown University, Providence, RI 02912 USA.
  ORCID: 0000-0003-3887-3248}

\thanks{
Pasten's research was supported by ANID FONDECYT Regular Grant 1230507 from Chile. Silverman's research was supported by Simons Collaboration Grant \#712332. Both authors' research was supported by the National Science Foundation under Grant No.\ 1440140 while the authors were in residence at the Mathematical Sciences Research Institute in Berkeley, California, during the spring of 2023. 
}


\newcommand{\JOE}[1]{{\textup{\color{blue} $\bigstar$ \textsf{{\bf Joe:} [#1]}}}}
\newcommand{\HECTOR}[1]{{\textup{\color{blue} $\bigstar$ \textsf{{\bf Hector:} [#1]}}}}

\hyphenation{ca-non-i-cal semi-abel-ian}


\newtheorem{theorem}{Theorem}[section]
\newtheorem{lemma}[theorem]{Lemma}
\newtheorem{sublemma}[theorem]{Sublemma}
\newtheorem{conjecture}[theorem]{Conjecture}
\newtheorem{protoconjecture}[theorem]{Proto-Conjecture}
\newtheorem{proposition}[theorem]{Proposition}
\newtheorem{corollary}[theorem]{Corollary}

\theoremstyle{definition}
\newtheorem*{claim}{Claim}
\newtheorem{definition}[theorem]{Definition}
\newtheorem*{intuition}{Intuition}
\newtheorem{example}[theorem]{Example}
\newtheorem{remark}[theorem]{Remark}
\newtheorem{question}[theorem]{Question}

\theoremstyle{remark}
\newtheorem*{acknowledgement}{Acknowledgements}
\newtheorem{exercise}[theorem]{exercise}

\newcommand{\TABT}[1]{\begin{tabular}[t]{@{}l@{}}#1\end{tabular}}
\newcommand{\TAB}[1]{\begin{tabular}{@{}l@{}}#1\end{tabular}}
\newcommand{\TABC}[1]{\begin{tabular}{@{}c@{}}#1\end{tabular}}

\newcommand{\Imply}[3]
{\label{imply:#1} \textup{(#2)}&\;\Longrightarrow\textup{(#3)}}
\newcommand{\NotImply}[3]
{\label{imply:#1} \textup{(#2)}&\;\centernot\Longrightarrow\textup{(#3)}}
\newcommand{\ImplyP}[3]{\par\noindent\eqref{imply:#1}\enspace
\framebox{$\textup{(#2)}\;\Longrightarrow\textup{(#3)}$}\enspace}
\newcommand{\NotImplyP}[3]{\par\noindent\eqref{imply:#1}\enspace
\framebox{$\textup{(#2)}\;\centernot\Longrightarrow\textup{(#3)}$}\enspace}


\newenvironment{notation}[0]{%
  \begin{list}%
    {}%
    {\setlength{\itemindent}{0pt}
     \setlength{\labelwidth}{4\parindent}
     \setlength{\labelsep}{\parindent}
     \setlength{\leftmargin}{5\parindent}
     \setlength{\itemsep}{0pt}
     }%
   }%
  {\end{list}}

\newenvironment{parts}[0]{%
  \begin{list}{}%
    {\setlength{\itemindent}{0pt}
     \setlength{\labelwidth}{1.5\parindent}
     \setlength{\labelsep}{.5\parindent}
     \setlength{\leftmargin}{2\parindent}
     \setlength{\itemsep}{0pt}
     }%
   }%
  {\end{list}}
\newcommand{\Part}[1]{\item[\upshape#1]}

\newcommand{\Case}[2]{\paragraph{\textbf{\boldmath Case #1: #2.}}\hfil\break\ignorespaces}
\newcommand{\Claim}[2]{\paragraph{\textbf{\boldmath Claim #1: #2}}\hfil\break\ignorespaces}

\renewcommand{\a}{\alpha}
\newcommand{\bfalpha}{{\boldsymbol{\alpha}}}
\renewcommand{\b}{\beta}
\newcommand{\bfbeta}{{\boldsymbol{\beta}}}
\newcommand{\g}{\gamma}
\renewcommand{\d}{\delta}
\newcommand{\e}{\epsilon}
\newcommand{\f}{\varphi}
\newcommand{\bfphi}{{\boldsymbol{\f}}}
\renewcommand{\l}{\lambda}
\renewcommand{\k}{\kappa}
\newcommand{\lhat}{\hat\lambda}
\newcommand{\m}{\mu}
\newcommand{\bfmu}{{\boldsymbol{\mu}}}
\renewcommand{\o}{\omega}
\newcommand{\bfpi}{{\boldsymbol{\pi}}}
\renewcommand{\r}{\rho}
\newcommand{\bfrho}{{\boldsymbol{\rho}}}
\newcommand{\rbar}{{\bar\rho}}
\newcommand{\s}{\sigma}
\newcommand{\sbar}{{\bar\sigma}}
\renewcommand{\t}{\tau}
\newcommand{\z}{\zeta}
\newcommand{\D}{\Delta}
\newcommand{\G}{\Gamma}
\newcommand{\F}{\Phi}
\renewcommand{\L}{\Lambda}

\newcommand{\ga}{{\mathfrak{a}}}
\newcommand{\gA}{{\mathfrak{A}}}
\newcommand{\gb}{{\mathfrak{b}}}
\newcommand{\gB}{{\mathfrak{B}}}
\newcommand{\gc}{{\mathfrak{c}}}
\newcommand{\gM}{{\mathfrak{M}}}
\newcommand{\gn}{{\mathfrak{n}}}
\newcommand{\gp}{{\mathfrak{p}}}
\newcommand{\gP}{{\mathfrak{P}}}
\newcommand{\gS}{{\mathfrak{S}}}
\newcommand{\gq}{{\mathfrak{q}}}

\newcommand{\Abar}{{\bar A}}
\newcommand{\Ebar}{{\bar E}}
\newcommand{\kbar}{{\bar k}}
\newcommand{\Kbar}{{\bar K}}
\newcommand{\Pbar}{{\bar P}}
\newcommand{\Sbar}{{\bar S}}
\newcommand{\Tbar}{{\bar T}}
\newcommand{\gbar}{{\bar\gamma}}
\newcommand{\lbar}{{\bar\lambda}}
\newcommand{\ybar}{{\bar y}}
\newcommand{\phibar}{{\bar\f}}
\newcommand{\nubar}{{\overline\nu}}

\newcommand{\Acal}{{\mathcal A}}
\newcommand{\Bcal}{{\mathcal B}}
\newcommand{\Ccal}{{\mathcal C}}
\newcommand{\Dcal}{{\mathcal D}}
\newcommand{\Ecal}{{\mathcal E}}
\newcommand{\Fcal}{{\mathcal F}}
\newcommand{\Gcal}{{\mathcal G}}
\newcommand{\Hcal}{{\mathcal H}}
\newcommand{\Ical}{{\mathcal I}}
\newcommand{\Jcal}{{\mathcal J}}
\newcommand{\Kcal}{{\mathcal K}}
\newcommand{\Lcal}{{\mathcal L}}
\newcommand{\Mcal}{{\mathcal M}}
\newcommand{\Ncal}{{\mathcal N}}
\newcommand{\Ocal}{{\mathcal O}}
\newcommand{\Pcal}{{\mathcal P}}
\newcommand{\Qcal}{{\mathcal Q}}
\newcommand{\Rcal}{{\mathcal R}}
\newcommand{\Scal}{{\mathcal S}}
\newcommand{\Tcal}{{\mathcal T}}
\newcommand{\Ucal}{{\mathcal U}}
\newcommand{\Vcal}{{\mathcal V}}
\newcommand{\Wcal}{{\mathcal W}}
\newcommand{\Wbar}{{\overline{\mathcal W}}}
\newcommand{\Xcal}{{\mathcal X}}
\newcommand{\Ycal}{{\mathcal Y}}
\newcommand{\Zcal}{{\mathcal Z}}

\renewcommand{\AA}{\mathbb{A}}
\newcommand{\BB}{\mathbb{B}}
\newcommand{\CC}{\mathbb{C}}
\newcommand{\FF}{\mathbb{F}}
\newcommand{\EE}{\mathbb{E}}
\newcommand{\GG}{\mathbb{G}}
\newcommand{\NN}{\mathbb{N}}
\newcommand{\PP}{\mathbb{P}}
\newcommand{\QQ}{\mathbb{Q}}
\newcommand{\RR}{\mathbb{R}}
\newcommand{\TT}{\mathbb{T}}
\newcommand{\ZZ}{\mathbb{Z}}

\newcommand{\bfa}{{\boldsymbol a}}
\newcommand{\bfb}{{\boldsymbol b}}
\newcommand{\bfc}{{\boldsymbol c}}
\newcommand{\bfd}{{\boldsymbol d}}
\newcommand{\bfe}{{\boldsymbol e}}
\newcommand{\ee}{{\boldsymbol{e}}} 
\newcommand{\bff}{{\boldsymbol f}}
\newcommand{\bfg}{{\boldsymbol g}}
\newcommand{\bfi}{{\boldsymbol i}}
\newcommand{\bfj}{{\boldsymbol j}}
\newcommand{\bfk}{{\boldsymbol k}}
\newcommand{\bfm}{{\boldsymbol m}}
\newcommand{\bfn}{{\boldsymbol n}}
\newcommand{\bfp}{{\boldsymbol p}}
\newcommand{\bfr}{{\boldsymbol r}}
\newcommand{\bfs}{{\boldsymbol s}}
\newcommand{\bft}{{\boldsymbol t}}
\newcommand{\bfu}{{\boldsymbol u}}
\newcommand{\bfv}{{\boldsymbol v}}
\newcommand{\bfw}{{\boldsymbol w}}
\newcommand{\bfx}{{\boldsymbol x}}
\newcommand{\bfy}{{\boldsymbol y}}
\newcommand{\bfz}{{\boldsymbol z}}
\newcommand{\bfA}{{\boldsymbol A}}
\newcommand{\bfF}{{\boldsymbol F}}
\newcommand{\bfB}{{\boldsymbol B}}
\newcommand{\bfD}{{\boldsymbol D}}
\newcommand{\bfG}{{\boldsymbol G}}
\newcommand{\bfI}{{\boldsymbol I}}
\newcommand{\bfM}{{\boldsymbol M}}
\newcommand{\bfP}{{\boldsymbol P}}
\newcommand{\bfQ}{{\boldsymbol Q}}
\newcommand{\bfT}{{\boldsymbol T}}
\newcommand{\bfU}{{\boldsymbol U}}
\newcommand{\bfX}{{\boldsymbol X}}
\newcommand{\bfY}{{\boldsymbol Y}}
\newcommand{\bfzero}{{\boldsymbol{0}}}
\newcommand{\bfone}{{\boldsymbol{1}}}

\newcommand{\arithord}{\mathsf{A}}
\newcommand{\Aut}{\operatorname{Aut}}
\newcommand{\Berk}{{\textup{Berk}}}
\newcommand{\Birat}{\operatorname{Birat}}
\newcommand{\characteristic}{\operatorname{char}}
\newcommand{\CircleNum}[1]{\raisebox{.5pt}{\textcircled{\raisebox{-.9pt} {\small#1}}}}
\newcommand{\Closure}{\operatorname{\textsf{Cl}}}
\newcommand{\codim}{\operatorname{codim}}
\newcommand{\Count}{\operatorname{\textsf{N}}}
\newcommand{\CountSet}{\operatorname{\mathfrak{N}}}
\newcommand{\Crit}{\operatorname{Crit}}
\newcommand{\crit}{{\textup{crit}}}
\newcommand{\critwt}{\operatorname{critwt}} 
\newcommand{\Cycle}{\operatorname{Cycles}}
\newcommand{\dense}{{\textup{dense}}}
\newcommand{\diag}{\operatorname{diag}}
\newcommand{\dimEnd}{{M}}  
\newcommand{\dimpadic}{\dim_{\textup{$p$-adic}}}
\newcommand{\dimKrull}{\dim_{\textup{Krull}}}
\newcommand{\Disc}{\operatorname{Disc}}
\newcommand{\Div}{\operatorname{Div}}
\newcommand{\Df}{{Df}}  
\newcommand{\Dom}{\operatorname{Dom}}
\newcommand{\dyn}{{\textup{dyn}}}
\newcommand{\End}{\operatorname{End}}
\newcommand{\PortEndPt}{{\textup{endpt}}} 
\newcommand{\END}{\smash[t]{\overline{\operatorname{End}}}\vphantom{E}}
\newcommand{\EndPoint}{E}  
\newcommand{\Expectation}{\operatornamewithlimits{\mathbb{E}}} 
\newcommand{\ExtOrbit}{\mathcal{EO}} 
\newcommand{\Fbar}{{\bar{F}}}
\newcommand{\fgeq}{\mathbin{\equiv_f}}
\newcommand{\fib}{{\textup{fib}}}
\newcommand{\Fix}{\operatorname{Fix}}
\newcommand{\Fiber}{\operatorname{Fiber}}
\newcommand{\FOD}{\operatorname{FOD}}
\newcommand{\FOM}{\operatorname{FOM}}
\newcommand{\Frame}{\operatorname{Fr}}
\newcommand{\Gal}{\operatorname{Gal}}
\newcommand{\genus}{\operatorname{genus}}
\newcommand{\GITQuot}{/\!/}
\newcommand{\GL}{\operatorname{GL}}
\newcommand{\Gp}{\Gcal}
\newcommand{\GR}{\operatorname{\mathcal{G\!R}}}
\newcommand{\grand}{{\textup{grand}}} 
\newcommand{\full}{\pm} 
\newcommand{\hhat}{{\hat h}}
\newcommand{\hplus}{h^{\scriptscriptstyle+}}
\newcommand{\Hom}{\operatorname{Hom}}
\newcommand{\Index}{\operatorname{Index}}
\newcommand{\Image}{\operatorname{Image}}
\newcommand{\Isog}{\operatorname{Isog}}
\newcommand{\Isom}{\operatorname{Isom}}
\newcommand{\Jac}{\operatorname{Jac}}
\newcommand{\Ker}{{\operatorname{Ker}}}
\newcommand{\Ksep}{K^{\text{sep}}}  
\newcommand{\Length}{\operatorname{Length}}
\newcommand{\Lift}{\operatorname{Lift}}
\newcommand{\limstar}{\lim\nolimits^*}
\newcommand{\limstarn}{\lim_{\hidewidth n\to\infty\hidewidth}{\!}^*{\,}}
\newcommand\LS[2]{{\genfrac{(}{)}{}{}{#1}{#2}}} 
\newcommand{\Mat}{\operatorname{Mat}}
\newcommand{\maxplus}{\operatornamewithlimits{\textup{max}^{\scriptscriptstyle+}}}
\newcommand{\MOD}[1]{~(\textup{mod}~#1)}
\newcommand{\Model}{\operatorname{Model}}
\newcommand{\Mor}{\operatorname{Mor}}
\newcommand{\Moduli}{\mathcal{M}}
\newcommand{\MODULI}{\overline{\mathcal{M}}}
\newcommand{\Mult}{\operatorname{\textup{\textsf{Mult}}}}
\newcommand{\Norm}{{\operatorname{\mathsf{N}}}}
\newcommand{\notdivide}{\nmid}
\newcommand{\normalsubgroup}{\triangleleft}
\newcommand{\NS}{\operatorname{NS}}
\newcommand{\onto}{\twoheadrightarrow}
\newcommand{\ord}{\operatorname{ord}}
\newcommand{\Orbit}{\mathcal{O}}
\newcommand{\Pcase}[3]{\par\noindent\framebox{$\boldsymbol{\Pcal_{#1,#2}}$}\enspace\ignorespaces}
\newcommand{\Per}{\operatorname{Per}}
\newcommand{\Perp}{\operatorname{Perp}}
\newcommand{\PrePer}{\operatorname{PrePer}}
\newcommand{\PGL}{\operatorname{PGL}}
\newcommand{\Pic}{\operatorname{Pic}}
\newcommand{\prim}{\textup{prim}}
\newcommand{\Prob}{\operatorname{Prob}}
\newcommand{\Proj}{\operatorname{Proj}}
\newcommand{\Qbar}{{\bar{\QQ}}}
\newcommand{\QR}[2]{\left( \dfrac{#1}{#2} \right) }
\newcommand{\rank}{\operatorname{rank}}
\newcommand{\Rat}{\operatorname{Rat}}
\newcommand{\Resultant}{\operatorname{Res}}
\newcommand{\Residue}{\operatorname{Residue}} 
\renewcommand{\setminus}{\smallsetminus}
\newcommand{\sgn}{\operatorname{sgn}}
\newcommand{\SL}{\operatorname{SL}}
\newcommand{\Sing}{\operatorname{Sing}} 
\newcommand{\Span}{\operatorname{Span}}
\newcommand{\Spec}{\operatorname{Spec}}
\renewcommand{\ss}{{\textup{ss}}}
\newcommand{\stab}{{\textup{stab}}}
\newcommand{\Stab}{\operatorname{Stab}}
\newcommand{\Support}{\operatorname{Supp}}
\newcommand{\Sym}{\operatorname{Sym}}  
\newcommand{\tors}{{\textup{tors}}}
\newcommand{\Trace}{\operatorname{Trace}}
\newcommand{\trianglebin}{\mathbin{\triangle}} 
\newcommand{\tr}{{\textup{tr}}} 
\newcommand{\UHP}{{\mathfrak{h}}}    
\newcommand{\val}{\operatorname{val}} 
\newcommand{\Var}{\operatornamewithlimits{Var}} 
\newcommand{\wt}{\operatorname{wt}} 
\newcommand{\<}{\langle}
\renewcommand{\>}{\rangle}

\newcommand{\pmodintext}[1]{~\textup{(mod}~#1\textup{)}}
\newcommand{\ds}{\displaystyle}
\newcommand{\longhookrightarrow}{\lhook\joinrel\longrightarrow}
\newcommand{\longonto}{\relbar\joinrel\twoheadrightarrow}
\newcommand{\SmallMatrix}[1]{%
  \left(\begin{smallmatrix} #1 \end{smallmatrix}\right)}


\begin{abstract}
Let $X/K$ be a smooth projective variety defined over a number field, and let $f:X\to{X}$ be a morphism defined over $K$. We formulate a number of statements of varying strengths asserting, roughly, that if there is at least one point $P_0\in{X(K)}$ whose $f$-orbit $\mathcal{O}_f(P_0):=\bigl\{f^n(P):n\in\mathbb{N}\bigr\}$ is Zariski dense, then after replacing~$K$ by a finite extension, there are many~$f$-orbits in~$X(K)$.  For example, a weak conclusion would be that $X(K)$ is not the union of finitely many (grand) $f$-orbits, while a strong conclusion would be that any set of representatives for the Zariski dense grand $f$-orbits in~$X(K)$ is itself Zariski dense. We prove statements of this sort for various classes of varieties and maps, including projective spaces, abelian varieties, and surfaces.
\end{abstract}

\maketitle

\tableofcontents

\section{Introduction}
\label{section:introduction}
Let~$K$ be a number field, and let~$X/K$ be a smooth projective variety, and let
\[
f:X\to{X}
\]
be an endomorphism of defined over~$K$. The theme of this article is that if there is even a single point~$P_0\in{X(K)}$ whose forward $f$-orbit
\[
\Orbit_f(P_0) := \bigl\{ f^n(P_0) : n \ge 0 \bigr\}
\]
is Zariski dense in~$X$, then after replacing~$K$ by a finite extension, the set~$X(K)$ contains lots of distinct large orbits whose points are widely distributed. There are many ways to turn this vague idea, which we call an \emph{orbit propagation principle}, into a precise statement.  We will describe four orbit propagation principles here, and we refer the reader to Section~\ref{section:notdefconj}, and especially Table~\ref{table:propogationstatements}, for many others. In order to state these principles, we begin with some useful definitions and notation.

\begin{definition}
The set of points with Zariski dense~$f$-orbit is denoted
\[
X_f^\dense := \bigl\{ P\in X : \overline{\Orbit_f(P)}=X \bigr\}.
\]
\end{definition}

\begin{definition}
The \emph{grand orbit} of a point~$P\in{X}$ is the set of points whose orbits eventually merge with the orbit of~$P$,
\[
\Orbit_f^\grand(P) 
:= \bigl\{ Q\in X : \Orbit_f(P)\cap\Orbit_f(Q)\ne\emptyset\bigr\}.
\]
We say that~$P$ and~$Q$ are \emph{grand $f$-orbit equivalent}, and we write \text{$P\fgeq{Q}$}, if
\[
\Orbit_f(P)\cap\Orbit_f(Q)\ne\emptyset,
\quad\text{or equivalently, if}\quad
\Orbit_f^\grand(P) = \Orbit_f^\grand(Q). 
\]
\end{definition}

We refer the reader to Lemma~\ref{lemma:orbitelemproperties}(a) for a proof that grand $f$-orbit equivalence is an equivalence relation.

\begin{conjecture}
\label{conjecture:weakstrongOPP}
Assume that~$X(K)$ has at least one Zariski dense~$f$-orbit. Then there is a finite extension~$K'/K$ such that\textup:
\begin{parts}
\Part{(a)}
\textup{(Orbit Propagation Principle (B1))}
For all~$P_1,\ldots,P_r\in{X(K')}$, the set 
\[
X(K')\setminus\bigl( \Orbit_f(P_1) \cup\cdots\cup \Orbit_f(P_r) \bigr)
\]
is Zariski dense in~$X$.
\Part{(b)}
\textup{(Orbit Propagation Principle (C1))}
For all~$P_1,\ldots,P_r\in{X(K')}$, the set 
\[
X(K')\setminus\bigl( \Orbit_f^\grand(P_1) \cup\cdots\cup \Orbit_f^\grand(P_r) \bigr)
\]
is Zariski dense in~$X$.
\Part{(b$'$)}
\textup{(Orbit Propagation Principle (C2$\exists$))}
The set~$X(K')$ contains a Zariski dense set of
coset representatives for~$X(K')/\fgeq$.
\Part{(c)}
\textup{(Orbit Propagation Principle (C3$\forall$))}
Every set of coset representatives in~$X(K')$ for the quotient 
\[
X_f^\dense(K')/\fgeq
\quad\text{is Zariski dense in $X$.}
\]
\end{parts}
\end{conjecture}

\begin{intuition}
Conjectures~\ref{conjecture:weakstrongOPP}(a) and~(b) say that a finite union of \textup(grand\textup) orbits can cover only a small portion of the $K'$-rational points. Conjectures~\ref{conjecture:weakstrongOPP}(b$'$) and~(c) say that taking one (any) $K'$-rational point from each (Zariski dense) grand orbit in~$X(K')$ always results in a Zariski dense set.
\end{intuition}

To what extent do propagation principles of varying strengths imply one another? There are some implications that are universally true and (mostly) easy to prove. For example, in the notation of Conjecture~\ref{conjecture:weakstrongOPP}, we always have
\[
\text{(C3$\forall$)}
\quad\Longrightarrow\quad
\text{(C2$\exists$)}
\quad\Longleftrightarrow\quad
\text{(C1)}
\quad\Longrightarrow\quad
\text{(B1)}.
\]
See Table~\ref{table:propagationimplications} for a more extensive diagram of universal implications and Appendix~\ref{appendix:implicationspropstatements} for proofs.
\par
The bulk of this article is devoted to proving non-trivial orbit propagation properties of varying strengths from the initial weak assumption that there is a single Zariski dense~$f$-orbit. We have no general statement, but we prove results  or various classes of varieties and maps. We state here some exemplary results. We refer the reader to Theorem~\ref{theorem:orbitpropagationprovencases} for a complete description of the results in this paper and to Sections~\ref{section:projectivespace}--\ref{section:surfaces} for the proofs. 

\begin{theorem}
\label{theorem:exampleresults}
\leavevmode
\begin{parts}
\Part{(a)}
Orbit Propagation Principle~\textup{(B1)} \textup(Conjecture~\textup{\ref{conjecture:weakstrongOPP}(a)}\textup) is true
for endomorphisms of smooth projective surfaces.
\Part{(b)}
The equivalent
Orbit Propagation Principles~{\textup{(C1) and (C2$\exists$)}} \textup(Conjectures~\textup{\ref{conjecture:weakstrongOPP}(b,b$'$)}\textup) are true for endomorphisms of projective space~$\PP^N$.
\Part{(c)}
Orbit Propagation Principle~\textup{(C3$\forall$)} \textup(Conjecture~\textup{\ref{conjecture:weakstrongOPP}(c)}\textup) is true for  linear endomorphisms of~$\PP^N$ and for endomorphisms of geometrically simple abelian varieties.
\end{parts}
\end{theorem}
\begin{proof}
(a)\enspace
See Theorem~\ref{theorem:surfaces}.
\par\noindent(b)\enspace
See Theorem~\ref{theorem:projectivespacedeg2}(a).
\par\noindent(c)\enspace
See Theorem~\ref{theorem:projectivespacedeg2}(c) for~$\PP^N$ and Theorem~\ref{theorem:simpleabelianvariety} for abelian varieties.
\end{proof}

\begin{remark}
We note that although we are able to prove our strongest Orbit Propagation Principle~\text{(C3$\forall$)} (Conjecture~\ref{conjecture:weakstrongOPP}(c)) in only a limited number of cases, we know of no examples for which it fails to be true.
\end{remark}

We now describe the structure of this paper. We start in Section~\ref{section:notdefconj} with definitions, notation, the description of a number of different orbit propagation principles, and a statement of our main results. Section~\ref{section:motivation} describes in more detail the motivation that led to the idea of orbit propagation. Then Sections~\ref{section:projectivespace}--\ref{section:surfaces} contain the proofs of our main results, using a variety of techniques and tools that include canonical heights and height counting functions, $p$-adic methods, algebro-geometric techniques, and deep theorems of Faltings et al.\ on the intersection of subvarieties of abelian varieties with subgroups of finite type. We include two appendices. 
Appendix~\ref{section:ZDconjecture} briefly discusses various versions of the related Zariski Density Conjecture and gives pointers to the literature. Appendix~\ref{appendix:elementaryresults} contains a number of auxiliary results, including elementary properties of orbits (Lemma~\ref{lemma:orbitelemproperties}) in Section~\ref{appendix:elempropsorbits}, elementary implications relating the various orbit propagation properties (Proposition~\ref{proposition:propagationimplications} and Lemma~\ref{lemma:C2forallZDimpliesC3forall}) in Section~\ref{appendix:implicationspropstatements}, and a weak height counting estimate for~$\PP^N(K)$ (Lemma~\ref{lemma:ctYKTleCTN}) in Section~\ref{appendix:htcountPN}.

We conclude this introduction with a few additional remarks.

\begin{figure}[ht]
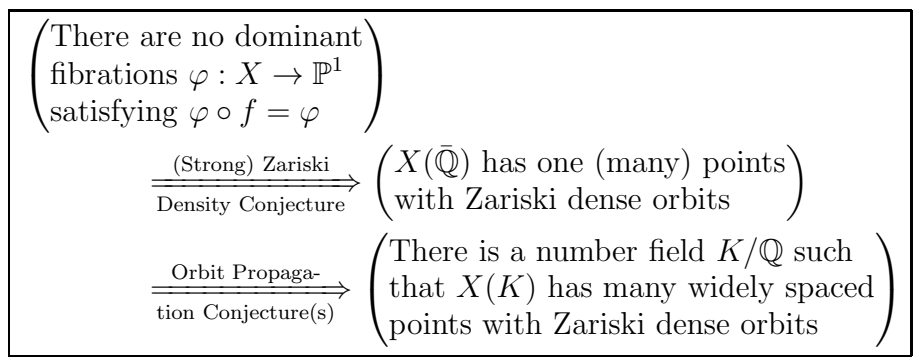

\framebox{$
\begin{aligned}
&\left(\text{\TAB{
There are no dominant\\
fibrations $\f:X\to\PP^1$\\
satisfying $\f\circ f=\f$\\
}}\right) 
\\
&\hspace*{4em}
\xRightarrow[\text{Density Conjecture}]{\text{(Strong) Zariski}}
\left(\text{\TAB{
$X(\Qbar)$ has one (many) points\\
with Zariski dense orbits\\
}}\right)
\\
&\hspace*{4em}
\xRightarrow[\text{tion Conjecture(s)}]{\text{Orbit Propaga-}}
\left(\text{\TAB{
There is a number field $K/\QQ$ such\\
that $X(K)$ has many widely spaced\\
points with Zariski dense orbits\\
}}\right)
\end{aligned}
$}
\caption{Relation between the Zariski Density Conjectures and the Orbit Propagation Conjectures}
\label{figure:ZDCimpliesOPC}
\end{figure}

\begin{remark}[The Zariski Density Conjecture]
There is a large literature on the following \emph{Zariski density conjecture} of Amerik--Bogomolov--Rovinsky and Medvedev--Scanlon (see also Amerik~\cite{MR2784670} and Zhang~\cite{MR2408228}). We briefly describe the relation of the Zariski density conjecture to the present paper, with further discussion and details in Section~\ref{section:ZDconjecture}.

\begin{conjecture}[The Zariski-Density Orbit Conjecture]
\label{conjecture:ZDconjecture1}
\textup{(\cite[Conjecture~1.2]{MR2862064} and \cite[Conjecture 5.10]{MR3126567})}  
Let~$X/\Qbar$ be a smooth variety, and let~$f:X\dashrightarrow{X}$ be a dominant rational map. Then one of the following is true\textup:
\begin{parts}
\Part{(a)}
There is a point~$P\in{X(\Qbar)}$ whose orbit~$\Orbit_f(P)$ is well-defined and Zariski dense in~$X$.
\Part{(b)}
There is a dominant rational map~$\f:X\dashrightarrow\PP^1$ such that~$\f\circ{f}=f$. \textup(Note that in this case, the orbits of~$f$ are restricted to lie in the fibers of~$\f$, and thus~$f$ cannot have any Zariski dense orbits.\textup)
\end{parts}
\end{conjecture}

Xie and others have formulated stronger versions of Conjecture~\ref{conjecture:ZDconjecture1}, and in particular, Xie~\cite[Section~1.3]{MR2862064} notes that ``if we have one Zariski dense orbit, we expect many such orbits.'' However, both the original conjecture and its various generalizations assert only that~$X(\Qbar)$ contains ``many'' points with Zariski dense orbit. They do not appear to consider the problem of whether there is a finite extension~$K/\QQ$ such that~$X(K)$ contains many such points; see Section~\ref{section:ZDconjecture} for details. Thus the Zariski Density Conjectures and our Orbit Propagation Conjectures are complementary. The implication relationship, at least for endomorphisms of non-singular projective varieties, is summarized in Figure~\ref{figure:ZDCimpliesOPC}.
\end{remark}

\begin{remark}[Motivation]
Our motivation for an orbit propagation principle arose from a~$35$-year old conjecture~\cite{MR1009803} of the first author. Very roughly, the conjecture says that if~$X(K)$ has infinitely many points, then ignoring ``error terms,'' the height counting function of~$X(K)$ should grow like a power of~$T$, a power of~$\log(T)$, or be bounded. (See Section~\ref{section:motivation} for a precise statement.) On the other hand, a conjecture of Kawaguchi and the second author~\cite{MR3456169} suggests that if~$f$ is sufficiently dynamically complicated (formally, if the dynamical degree of~$f$ satisfies $\d(f)>1$), then the height counting function of~$\Orbit_f(P)$ of the orbit of a point grows no faster than~$\log\log(T)$. Hence~$X(K)$ should have lots of different orbits. This vague idea led first to various weak conjectures such as~(A) implies~(B1) and~(C1) in Table~\ref{table:propogationstatements}, and eventually to stronger statements, culminating in the strong orbit propagation principle described in Conjecture~\ref{conjecture:weakstrongOPP}(c).
\end{remark}

\begin{remark}[Rational Self-Maps of Singular Varieties?]
\label{remark:ratmapsingularvar}
In this paper we consider orbit propagation principles only for self-morphisms~$f$ of smooth projective varieties~$X$. One might ask whether similar statements hold for rational maps and/or for non-smooth quasi-projective varieties, as is allowed for example in the Zariski density conjecture (Conjecture~\ref{conjecture:ZDconjecture1}).
We leave such questions for the future.
\end{remark}

\begin{acknowledgement}
We would like to thank Brendan Hassett, Yohsuke Matsuzawa, Max Weinreich, Junyi Xie, and De-Qi Zhang for helpful comments and assistance. 
\end{acknowledgement}

\section{Notation, Definitions, Conjectures and Main Results}
\label{section:notdefconj}

\begin{definition}
\label{definition:orbits}
Throughout this article, we fix the following notation:
\begin{align*}
    X&\quad\text{a smooth projective variety with $\dim(X)\ge1$.} \\
    f&\quad\text{an endomorphism $f:X\to X$.}
\end{align*}
Let~$P\in{X}$ be a geometric point of~$X$. We define various sorts of orbits of~$P$.\footnote{Our set of natural numbers~$\NN$ includes~$0$, so~$\Orbit_f(P)$ includes the point~$P=f^0(P)$.}
\begin{description}
\item[(Forward) $f$-orbit]
\[
\Orbit_f(P) = \Orbit^+_f(P) := \bigl\{ f^n(P) : n\in\NN \bigr\}.
\]
\item[Backward $f$-orbit]
\[
\Orbit^-_f(P) := \bigl\{ Q\in X : P\in\Orbit_f(Q) \bigr\}.
\]
\item[Full $f$-orbit]
\[
\Orbit^\full_f(P) := \Orbit_f^+(P) \cup \Orbit_f^-(P).
\]
\item[Grand $f$-orbit]
\[
\Orbit_f^\grand(P) = \bigl\{ Q\in X : \Orbit_f(P)\cap\Orbit_f(Q)\ne\emptyset \bigr\}.
\]
\end{description}
\end{definition}

\begin{definition}
Let~$P,Q\in{X}$. If 
\[
\Orbit_f(P)\cap\Orbit_f(Q)\ne\emptyset,
\]
then we say that~$P$ and~$Q$ are
\[
\text{\emph{grand $f$-orbit equivalent},\quad and we write\quad $P\fgeq{Q}$.}
\]
\end{definition}

\begin{definition}
We denote the set of points with Zariski dense $f$-orbit by
\[
X_f^\dense := \bigl\{ P\in X : \overline{\Orbit_f(P)}=X \bigr\}.
\]
\end{definition}

We refer the reader to Section~\ref{appendix:elementaryresults} for a proof that 
grand $f$-orbit equivalence is an equivalence relation on the set of points of~$X$, and for various other elementary properties of the different types of orbits.

\begin{conjecture}[\textup{Orbit Propagation}]
\label{conjecture:propagation}
Let~$\Xcal$ be a set of \textup(smooth projective connected\textup) varieties defined over~$\Qbar$, and for each~$X\in\Xcal$, let~$\Fcal_X$ be a collection of~$\Qbar$-morphisms~$X\to{X}$. We say that~$(\Xcal,\Fcal)$ has an \emph{orbit propagation property} if for every~$X\in\Xcal$ and every~$f\in\Fcal_X$, the following implication holds\textup:
\[
\left(\begin{tabular}{@{}l@{}}
$X(\Qbar)$ contains at\\
least one Zariski\\
dense $f$-orbit\\
\end{tabular}\right)
\quad\Longrightarrow\quad
\left(
\text{\parbox{.5\hsize}{\raggedright
there is a number field $K/\QQ$
that is a field of definition
for $X$ and $f$ such that $X(K)$
contains ``many large'' $f$-orbits
}}\right).
\]
The words ``many large'' may be quantified using the various orbit statements in Table~\textup{\ref{table:propogationstatements}}. So we say that~$(\Xcal,\Fcal)$ has an  \emph{orbit propagation property of a specified type} if for every~$X\in\Xcal$ and every~$f\in\Fcal_X$, it satisfies\textup:
\[
X_f^\dense(\Qbar)\ne\emptyset
\quad\Longrightarrow\quad
\left(
\text{\parbox{.6\hsize}{\raggedright
there is a number field~$K/\QQ$ that is a field of definition for~$X$ and~$f$ such that a specified~\textup{(B)} or~\textup{(C)} statement in Table~\textup{\ref{table:propogationstatements}} is valid for~$X(K)$
}}
\right).
\]
\end{conjecture}

\begin{remark}
We note that since we assume that~$X$ is a smooth projective variety and that~$f:X\to{X}$ is a morphism, assumption~(A) that there exists at least one Zariski dense orbit forces~$f$ to be surjective, and hence finite.
\end{remark}

\begin{remark}
\label{remark:(A)implies}
When we refer to Table~\ref{table:propogationstatements} and make an assertion such as
\[
\text{(A)} \quad\Longrightarrow\quad \text{(C1)},
\]
what we always mean is that if there is a number field~$K$ such that~(A) is true, then possibly after replacing~$K$ with a finite extension, the statement~(C1) is also true.
\end{remark}

\begin{table}[t]
\framebox{\vbox{%
\begin{parts}
\Part{}
{\centering\textbf{(A)-Statement: One Forward Orbit}\par}
\Part{(A)\enspace}
There is at least one Zariski dense $f$-orbit in
$X(K)$, i.e., $$X_f^\dense(K)\ne\emptyset.$$
\par\vspace{-7pt}\hbox to.98\hsize{\hrulefill}
\Part{}
{\centering\textbf{(B)-Statements: Many Forward Orbits}\par}
\Part{(B1)\enspace}
For any finite collection of~$f$-orbits~$\G_1,\ldots,\G_r$, the set
\[
X(K) \setminus (\G_1\cup\cdots\cup\G_r)
\quad\text{is Zariski dense in $X$.}
\]
\Part{(B1$\infty$)\enspace}
For every proper Zariski closed set $Y\subsetneq{X}$, the set
\[
X(K) \setminus \bigcup_{P\in Y(K)} \Orbit_f(P)
\quad\text{is Zariski dense in $X$.}
\]
\par\vspace{0pt}\hbox to.98\hsize{\hrulefill}
\Part{}
{\centering\textbf{(C)-Statements: Many Grand Forward Orbits}\par}
\Part{(C1)\enspace}
For any finite collection of grand~$f$-orbits~$\G_1,\ldots,\G_r$, the set
\[
X(K) \setminus (\G_1\cup\cdots\cup\G_r)
\quad\text{is Zariski dense in $X$.}
\]
\Part{(C1$\infty$)\enspace}
For every proper Zariski closed set $Y\subsetneq{X}$, the set
\[
X(K) \setminus \bigcup_{P\in Y(K)} \Orbit_f^\grand(P)
\quad\text{is Zariski dense in $X$.}
\]
\Part{(C2$\exists$)\enspace}
There exists a Zariski dense set of representatives in $X(K)$ for 
\[
X(K)/\fgeq.
\]
\Part{(C2$\forall$)\enspace}
Every complete set of representatives for 
\[
\text{$X(K)/\fgeq$\quad is Zariski dense in $X$.}
\]
\Part{(C3$\exists$)\enspace}
There exists a Zariski dense set of representatives in $X_f^\dense(K)$ for
\[
X_f^\dense(K)/\fgeq.
\]
\Part{(C3$\forall$)\enspace}
$X_f^\dense(K)\ne\emptyset$ and every complete set of representatives for 
\[
\text{$X_f^\dense(K)/\fgeq$\quad is Zariski dense in $X$.}
\]
\end{parts}
}}
\caption{Orbit propagation statements}
\label{table:propogationstatements}
\end{table}

The many orbit propagation statements in Table~\ref{table:propogationstatements} are related by a number of straightforward implications,\footnote{Although to be strictly accurate, the fact that~(C1) implies~(C2$\exists$) is not completely straightforward. In particular, it depends on the fact that our fields are countable.} which we describe pictorially in Table~\ref{table:propagationimplications} and prove in Proposition~\ref{proposition:propagationimplications}. 
In particular, we note that
\[
\text{(C3$\forall$)}
=
\text{``The one orbit propagation property to rule them all!''}
\]
So one might ask whether 
\begin{equation}
\label{eqn:rulethemall}
\text{(A)}\quad
\stackrel{{?}{?}{?}}{\Longrightarrow}
\quad\text{(C3$\forall$)}.
\end{equation}
We do not know of any examples for which~\eqref{eqn:rulethemall} fails to be true, and although we only able to prove it in certain cases, there are many classes of varieties and maps for which we can prove weaker orbit propagation implications. Theorem~\ref{theorem:orbitpropagationprovencases} summarizes our results.

\begin{table}
\framebox{\vbox{\hbox{%
$
\def\P#1#2{\textup{(#1$#2$)}}
\def\L{\Longrightarrow}
\def\D{\Big\Downarrow}
\begin{array}{*7c}
\P{C3}{\forall} & \L & \P{C2}{\forall} 
& \Longleftrightarrow & \P{C1}{\infty} & \L & \P{B1}{\infty} \\[1.5\jot]
\D && \D && \D && \D \\[2.5\jot]
\P{C3}{\exists} & \L & \P{C2}{\exists} & \Longleftrightarrow & \P{C1}{} & \L & \P{B1}{} \\[1.5\jot]
\D \\[2.5\jot]
\P{A}{} \\
\end{array}
$
}}}
\caption{\strut Universal implications relating the orbit propagation properties in Table~\ref{table:propogationstatements}.  See Proposition~\ref{proposition:propagationimplications}.}
\label{table:propagationimplications}
\end{table}

\begin{theorem}
\label{theorem:orbitpropagationprovencases}
The following orbit propagation implications hold for the indicated varieties and maps, with the proviso from Remark~\ref{remark:(A)implies} that one may need  to replace the field~$K$ with a finite extension\textup:
\begin{parts}
    \Part{(1)}
    Projective space, $\deg(f)\ge2$\textup: \\
    \textup{{(A) $\Rightarrow$ (C1) $\Leftrightarrow$ (C2$\exists$)}} and \textup{{(A) $\Rightarrow$ (B1$\infty$)}}  ---
    \textup{Theorem \ref{theorem:projectivespacedeg2}(a,b)}
    \Part{(2)}
    Projective space, $\deg(f)=1$\textup: \\
    \textup{{(A) $\Rightarrow$ (C3$\forall$)}} ---
    \textup{Theorem \ref{theorem:projectivespacedeg2}(c)}
    \Part{(3)}
    K3 Surfaces\textup: \\
    \textup{{(A) $\Rightarrow$ (C1) $\Leftrightarrow$ (C2$\exists$)}} ---
    \textup{Theorem \ref{theorem:K3}}
    \Part{(4)}
    Geometrically simple abelian varieties\textup: \\
    \textup{{(A) $\Rightarrow$ (C3$\forall$)}} ---
    \textup{Theorem \ref{theorem:simpleabelianvariety}}
    \Part{(5)}
    Smooth rational varieties and {\'e}tale $f$\textup: \\
    \textup{{(A) $\Rightarrow$ (B1)}} ---
    \textup{Corollary~\ref{corollary:CoroRat}}
    \Part{(6)}
    {\'E}tale quotients of abelian varieties\textup: \\
    \textup{{(A) $\Rightarrow$ (B1)}} ---
    \textup{Corollary~\ref{corollary:B1foretaleabelianquotients}}
    \Part{(7)}
    Smooth projective surfaces\textup: \\
    \textup{{(A) $\Rightarrow$ (B1)}} ---
    \textup{Theorem~\ref{theorem:surfaces}}
 \end{parts}
\end{theorem}

\section{Motivation for the Propagation Conjecture}
\label{section:motivation}
Our motivation for formulating some sort of propagation conjecture rests on an older, highly speculative, conjecture of the second author. That conjecture says roughly that, up to lower order terms, the height counting function for the integral points on an algebraic variety can have only one of three possible growth rates. The precise formulation requires some care balancing extending the field and discarding Zariski closed sets having fast growth rates. We give a precise statement for projective varieties and $K$-rational points; see the cited reference for a general formulation for quasi-projective varieties and~$S$-integral points.

\begin{definition}
\label{definition:arithmeticorder}
Let~$X/K$ be a smooth projective variety, and let
\[
H:X(\Kbar)\to[1,\infty)
\]
be a Weil height function associated to an ample divisor. We say that~$X$ has \emph{arithmetic order~$\arithord(X)$} if there is an integer~$m\ge2$ and a non-empty Zariski open subset~$U_0\subseteq{X}$ such that for every non-empty Zariski open subset~$U\subseteq{U_0}\subseteq{X}$ there is a finite extension~$L_0/K$ such that for every finite extension~$L/L_0$,
\begin{equation}
\label{eqn:arithorderlimit}
\lim_{T\to\infty} \frac{ \log^{(m)} \#\bigl\{ P\in U(L) : H(P) \le T\bigr\} }
{ \log^{(m+\arithord(X))} T } = 1.
\end{equation}
The notation~$\log^{(m)}$ denotes the~$m$-fold iterate, and by convention we set~$\arithord(X)=\infty$ if the limit is~$0$ for all~$U$ and all~$L$.
\end{definition}

Since the chain of logic in Definition~\ref{definition:arithmeticorder} is somewhat complicated, we note that it may be written succinctly using logical notation as
\[
\exists m\ge2,\;
\exists U_0\subseteq X,\;
\forall U\subseteq U_0,\;
\exists L_0/K,\;
\forall L/L_0,\;
\text{\eqref{eqn:arithorderlimit} is true.}
\]

\begin{conjecture}
\label{conjecture:Ulm}
\textup{\cite[Silverman 1987]{MR1009803}}
Let~$X/K$ be a smooth projective variety defined over a number field. Then the arithmetic order~$\arithord(X)$ exists and satisfies
\[
\arithord(X)=0 \quad\text{or}\quad \arithord(X)=1 \quad\text{or}\quad \arithord(X)=\infty.
\]
\end{conjecture}

\begin{remark}
We note that it is tempting to simply set~$m=2$ in Definition~\ref{eqn:arithorderlimit}, so we refer the reader to~\cite{MR1009803} for an example that suggests why it may be necessary in some cases to take~$m\ge3$. In any case, the conclusion of Conjecture~\ref{conjecture:Ulm} says roughly that one of the following is true, where we are being very coarse about ignoring error terms:
\begin{equation}
\label{eqn:growthtrichotomy}
\#\bigl\{ P\in X(K) : H(P) \le T \bigr\}
\quad
\begin{cases}
\text{grows like a power of $T$,} \\
\text{grows like a power of $\log T$,} \\
\text{is bounded at $T\to\infty$.} \\
\end{cases}
\end{equation}
\end{remark}

We next observe that in many situations, if $f:X\to{X}$ is an endomorphism of a smooth projective variety defined over a number field~$K$, and if~$P\in{X(P)}$, then the logarithmic height of the points in $f$-orbit of~$P$ tend to grow exponentially. Indeed, if the dynamical degree of~$f$ satisfies~$\d(f)>1$, and if~$\Orbit_f(P)$ is Zariski dense, then it is conjectured~\cite{MR3456169} that
\[
\lim_{n\to\infty} \sqrt[n]{\log H\bigl(f^n(P)\bigr)} = \d(f).
\]
If this is true, then the height counting function for the points in the orbit satisfy
\begin{equation}
\label{eqn:QHQTllloglogT}
\#\bigl\{ Q\in \Orbit_f(P) : H(Q) \le T \bigr\} \ll \log\log T,
\end{equation}
where the implied constant depends on~$f$ and~$P$, but is independent of~$T$. Thus if~$X(K)$ were to be the union of a finite number of $f$-orbits, at least one of which is Zariski dense, then its height counting function would not be bounded, yet would increase too slowly to satisfy the other growth conditions in Conjecture~\ref{conjecture:Ulm}. (This can be seen more clearly, albeit less precisely, by comparing~\eqref{eqn:growthtrichotomy} and~\eqref{eqn:QHQTllloglogT}.) Hence Conjecture~\ref{conjecture:Ulm} suggests that if~$X(K)$ contains at least one Zariski dense $f$-orbit, then (possibly after extending~$K$), it must contain infinitely many non-overlapping $f$-orbits. 
\par
We acknowledge that making a new conjecture on the basis of an older, not widely known, conjecture is somewhat dubious. Further,  using Conjecture~\ref{conjecture:Ulm} as a starting point, a natural conjecture would be a relatively weak statement such as the following:

\begin{protoconjecture}
Let~$K$ be a number field, let~$X/K$ be a smooth projective variety, let $f:X\to{X}$ be an endomorphism defined over~$K$ with dynamical degree~$\d(f)>1$. If~$\#X(K)=\infty$, then~$X(K)$ is not the union of finitely many~$f$-orbits.
\end{protoconjecture}

But as we explored this proto-conjecture, we realized that we were unable to rule out even much stronger statements, including dropping the~$\d(f)>1$ requirement, looking at grand orbits, requiring orbits to be Zariski dense, and changing the ``not a union of finitely many~$f$-orbits'' to a statement that complete sets of representatives for the~$f$-orbits must be Zariski dense. This led us to the plethora of orbit propagation properties listed in Table~\ref{table:propogationstatements}. The remainder of this paper is devoted to proving, for certain classes of varieties and maps, various versions of the statement that one Zariski dense $f$-orbit leads to many such orbits. 

\begin{remark}
We thank De-Qi Zhang for the following remarks and questions. 
\begin{parts}
\Part{(a)}
We have assumed that~$X$ is a smooth projective variety, but Conjecture~\ref{conjecture:weakstrongOPP} makes sense even if~$X$ is singular. 
\Part{(b)}
Suppose that~$f:X\to{X}$ descends or lifts to a surjective endomorphism $Y\to{Y}$ of a projective variety of the same dimension via a generically finite rational map or finite morphism $X\to{Y}$ or $Y\to{X}$. Is it true that Conjecture~\ref{conjecture:weakstrongOPP} is true for~$X$ if and only if it is true for~$Y$? If so, then the case of singular surfaces with~$\deg(f)\ge2$ can mostly be reduced to the case of smooth surfaces except possibly in some cases of polarized~$f$.
\Part{(c)}
To what extent can one prove Conjecture~\ref{conjecture:weakstrongOPP} or various other implications among the orbit propagation statements in Table~\ref{table:propogationstatements} if one assumes that~$f:X\to{X}$ is a polarized endomorphism, i.e., if there exists an ample line bundle~$\Lcal\in\Pic(X)\otimes\RR$ and a real number~$d>1$ such that~$f^*\Lcal\cong\Lcal^{\otimes{d}}$? We remark that examples of such morphisms include maps of degree~$\deg(f)\ge2$ on~$\PP^N$ and infinite order self-isogenies of abelian varieties, for which we have proven orbit propagation results in Theorems~\ref{theorem:projectivespacedeg2}(a,b) and~\ref{theorem:simpleabelianvariety}.
\end{parts}
\end{remark}
 
\section{Projective Space}
\label{section:projectivespace}

\begin{theorem}
\label{theorem:projectivespacedeg2}
Let $N\ge1$, let $f:\PP^N\to\PP^N$ be an endomorphism defined over~$\Qbar$, and assume that there is a point~$P_0\in\PP^N(\Qbar)$ whose orbit~$\Orbit_f(P_0)$ is Zariski dense in~$\PP^N$. There exists a number field~$K$ that is a common field of definition for~$f$ and~$P_0$ so that the following hold\textup:
\begin{parts}
\Part{(a)}
\framebox{$\deg(f)\ge2$}
For every finite collection of grand $f$-orbits 
\[
\G_1,\ldots,\G_r\subset\PP^N(\Qbar),
\]
the set
\[
{\PP^N(K) \setminus (\G_1\cup\cdots\cup\G_r)} 
\quad\text{is Zariski dense in $\PP^N$.}
\]
Equivalently, from Proposition~\ref{proposition:propagationimplications}, there exists a Zariski dense set in~$\PP^N(K)$ that contains one point in each grand $f$-orbit generated by the points in~$\PP^N(K)$. In the terminology of Table~\textup{\ref{table:propogationstatements}}, non-linear endomorphisms of~$\PP^N$ satisfy the orbit propagation statement
\[
\textup{(A)} \quad\Longrightarrow\quad \textup{(C1)}\;\Leftrightarrow\;\textup{(C2$\exists$)}.
\]
\Part{(b)}
\framebox{$\deg(f)\ge2$}
For every proper Zariski closed set $Y\subsetneq{\PP^N}$, the set
\[
\PP^N(K) \setminus \bigcup_{P\in Y(K)} \Orbit_f(P)
\quad\text{is Zariski dense in $\PP^N$.}
\]
In the terminology of Table~\textup{\ref{table:propogationstatements}}, non-linear endomorphisms of~$\PP^N$ satisfy the orbit propagation statement
\[
\textup{(A)} \quad\Longrightarrow\quad \textup{(B1$\infty$)}.
\]
\Part{(c)}
\framebox{$\deg(f)=1$}
Every complete set of representatives in $(\PP^N)_f^\dense(K)$ for $(\PP^N)_f^\dense(K)/\fgeq$ is Zariski dense in~$\PP^N$. 
In the terminology of Table~\textup{\ref{table:propogationstatements}}, linear endomorphisms of~$\PP^N$ satisfy the orbit propagation statement
\[
\textup{(A)} \quad\Longrightarrow\quad \textup{(C3$\forall$)}.
\]
\end{parts}
\end{theorem}

\begin{proof}
We remark that the proof of Theorem~\ref{theorem:projectivespacedeg2} uses three lemmas whose statements and proofs we defer until the end of this section.
\par\noindent(a,b)\enspace\framebox{$\deg(f)\ge2$}\enspace
The assumption that $\deg(f)\ge2$ means that there is a canonical height function~\cite{MR1255693}
\[
\hhat_f:\PP^N(\Kbar)\longrightarrow\RR_{\ge0}
\]
satisfying
\begin{align}
\label{eqn:canht1}
\hhat_f\circ f&=d\cdot\hhat_f,\\
\label{eqn:canht2}
|\hhat_f - h| &\le \Cl{hfh}(f),\\
\label{eqn:canht3}
\hhat_f(P)=0&\quad\Longleftrightarrow\quad \#\Orbit_f(P)<\infty.
\end{align}
\par
Combining~\eqref{eqn:canht2} and~\eqref{eqn:canht3} with the fact that~$\PP^N(K)$ has only finitely many points of bounded height, we see that~$\hhat_f$ takes on a minimal positive value, which we denote by
\begin{equation}
\label{eqn:hhatfminPNK}
\hhat_f^{\min}(\PP^N,K) 
:= \inf_{\substack{P\in\PP^N(K)\\\#\Orbit_f(P)=\infty\\}} \hhat_f(P)
> 0. 
\end{equation}
\par
Let~$\G$ be an $f$-orbit or a grand $f$-orbit such that~$\G\cap\PP^N(K)$ is non-empty. We define
\[
\hhat_f^{\min}(\G,K) = \inf\bigl\{ \hhat_f(Q) : Q\in\G\cap\PP^N(K) \bigr\}.
\]
Using~\eqref{eqn:canht2} and the fact that there are only finitely many points in~$\PP^N(K)$ of bounded Weil height, we see that there exists a point
\[
Q_{\G,K}\in\G\cap\PP^N(K)
\]
(not necessarily unique if~$\G$ is a grand orbit) such that
\[
\hhat_f(Q_{\G,K}) = \hhat_f^{\min}(\G,K).
\]
Further, we see that 
\begin{equation}
\label{eqn:hhatfminGKgt0}
\hhat_f^{\min}(\G,K)>0
\quad\Longleftrightarrow\quad
\text{$\G$ contains an $f$-wandering point}.
\end{equation}
(This is equivalent to every point in~$\G$ being $f$-wandering.)
\par\noindent(a)\enspace
Let~$\G$ be a grand $f$-orbit. Then
\begin{align*}
Q & \in\G\cap\PP^N(K) \\
&\quad\Longrightarrow\quad
f^i(Q)=f^j(Q_{\G,K})\quad\text{for some $i,j\in\NN$,} \\
&\quad\Longrightarrow\quad
d^i\cdot\hhat_f(Q)=d^j\cdot\hhat_f(Q_{\G,K})
\quad\text{for some $i,j\in\NN$,} \\
&\quad\Longrightarrow\quad
d^i\cdot\hhat_f(Q)=d^j\cdot\hhat_f(Q_{\G,K})
\quad\begin{tabular}[t]{@{}l@{}}
with $j\ge i$, 
since $\hhat_f(Q_{\G,K})$\\
is the smallest
value of~$\hhat_f$\\ for the points in~$\G$,\\
\end{tabular} \\
&\quad\Longrightarrow\quad
\hhat_f(Q) \in d^\NN \cdot\hhat_f^{\min}(\G,K)
\quad \text{since $\hhat_f(Q_{\G,K})=\hhat_f^{\min}(\G,K)$,} \\
&\quad\Longrightarrow\quad
\smash[b]{ h(Q) \in \bigcup_{n\in\NN} }
\Bigl[d^n\hhat_f^{\min}(\G,K)-\Cr{hfh}(f),d^n\hhat_f^{\min}(\G,K)+\Cr{hfh}(f)\Bigr] \\*
&\omit\hfill\quad\text{since $|\hhat_f-h|\le\Cr{hfh}(f)$.}
\end{align*}
\par
We now suppose that~$\G_1,\ldots,\G_r$ are grand orbits, and to ease notation, we let
\[
\Cl{G1Gr}(i) := \Cr{G1Gr}(f,\G_i,K) = \hhat_f^{\min}(\G_i,K).
\]
The above calculation shows that
\[
\bigcup_{i=1}^r \G_i\cap\PP^N(K)
\subseteq
\bigcup_{i=1}^r \bigcup_{n\in\NN}
\Bigl\{ Q\in\PP^N(K) : \bigl| h(Q) - d^n\Cr{G1Gr}(i) \bigr| \le \Cr{hfh}(f) \Bigr\}.
\]
Hence taking heights, we find that the set of heights of the points in the union of the \text{$\G_i\cap\PP^N(K)$} is contained in a union of intervals,
\begin{multline*}
\bigcup_{i=1}^r 
\Bigl\{ h(Q) : Q\in\G_i\cap\PP^N(K) \Bigr\}  \\
\subseteq \bigcup_{i=1}^r \bigcup_{n\in\NN}
\Bigl[ d^n\Cr{G1Gr}(i) - \Cr{hfh} ,\; d^n\Cr{G1Gr}(i) + \Cr{hfh} \Bigr].
\end{multline*}
An elementary estimate shows that the double union on the right-hand side omits intervals in~$\RR_{\ge0}$ of arbitrarily large length; see Lemma~\ref{lemma:missedinterval} for a more precise result. In particular, we can find infinitely many intervals of length~$1$ that are omitted, say
\begin{multline*}
0\le t_1 < t_2 < t_3 < \cdots\quad\text{satisfying}\quad
t_{i+1}>t_i+1\quad\text{and}\\
\biggl( \bigcup_{i\ge1} [t_i,t_i+1] \biggr) \;\cap\;
\biggl( \bigcup_{i=1}^r 
\Bigl\{ h(Q) : Q\in\G_i\cap\PP^N(K) \Bigr\} \biggr) = \emptyset.
\end{multline*}
Every interval of length~$1$ contains a number of the form~$\log(a)$ with~$a\in\NN$, so we can find a sequence of distinct positive integers~$a_1,a_2,\ldots$ satisfying
\[
\log(a_j) \notin \bigcup_{i=1}^r \Bigl\{ h(Q) : Q\in\G_i\cap\PP^N(K) \Bigr\}
\quad\text{for all $j\ge1$.}
\]
We consider the set of points
\[
\Acal := \bigcup_{j\ge1} \left\{ [a_j,b_1,\cdots,b_N]\in\PP^N(\QQ) : 
\begin{array}{@{}l@{}}
b_1,\ldots,b_N\in\ZZ,\\
0\le |b_1|,\ldots,|b_N|\le a_j\\
\end{array}\right\} .
\]
The heights of the points in~$\Acal$ are all of the form~$\log(a_j)$, so they are not heights of~$K$-rational points in any of the grand orbits~$\G_1,\ldots,\G_r$. This proves that
\[
\Acal\; \cap \; \bigcup_{i=1}^r \Bigl(\G_i\cap\PP^N(K) \Bigr)
= \emptyset.
\]
On the other hand, it is clear that~$\Acal$ is Zariski dense in~$\PP^N$. Hence $\PP^N(K)$ contains a Zariski dense set of points not lying in any of the grand orbits~$\G_1,\ldots,\G_r$, which completes the proof of orbit propagation property~(C1).
\par\noindent(b)\enspace
To simplify formulas, we are going to use Weil and canonical heights relative to the field~$K$. For any subset~$\Pcal\subseteq\PP^N(K)$, we define a counting function
\begin{equation}
\label{eqn:countingfunction}
\Count(\Pcal,T) := \#\bigl\{P\in\Pcal : H(P) \le T\bigr\}.
\end{equation}
With our height normalization, Lemma~\ref{lemma:ctYKTleCTN}(a) tells us that\footnote{Schanuel~\cite{MR557080} gives a formula, with error term, for~$\Count\bigl(\PP^N(K),T\bigr)$, but we will not need anything that precise.}
\begin{equation}
\label{eqn:countPNKT}
\Cl{PNKT}(K,N) T^{N+1} \le \Count\bigl(\PP^N(K),T\bigr) \le \Cl{PNKTup}(K,N) T^{N+1}.
\end{equation}
\par
Let~$P\in{X(K)}$ be an $f$-wandering point. Then
\begin{align}
\Count\bigl(\Orbit_f(P),T\bigr)
&= \#\bigl\{ n\ge0 : H\bigl(f^n(P)\bigr) \le T \bigr\} \notag \\
&\le \#\Bigl\{ n\ge0 : \hhat_f\bigl(f^n(P)\bigr) \le \log(T)+\Cr{hfh}(f) \Bigr\}
\quad\text{from \eqref{eqn:canht2},} \notag \\
&\le \#\Bigl\{ n\ge0 : d^n\hhat_f(P) \le \log(T)+\Cr{hfh}(f) \Bigr\} \notag \\
&\le 1 + \log_d \left( \frac{\log T + \Cr{hfh}(f)}{\hhat_f(P)} \right) \notag \\
&\le 1 + \log_d \left( \frac{\log T + \Cr{hfh}(f)}{\hhat_f^{\min}(\PP^N,K)} \right) 
\quad\text{from \eqref{eqn:hhatfminPNK},} \notag \\
\label{eqn:ctfororb2}
&\le \Cl{Oplus}(K,N,f)\cdot\log\log(T).
\end{align}
\par
We also note that
\begin{align*}
h(P) & > \log(T) + 2\Cr{hfh}(f)\\
&\quad\Longrightarrow\quad
\hhat_f(P) > \log(T) + \Cr{hfh}(f)
\quad\text{from \eqref{eqn:canht2},} \\
&\quad\Longrightarrow\quad
d^n\hhat_f(P) > d^n\log(T) + d^n\Cr{hfh}(f)
\quad\text{for all $n\in\NN$,} \\
&\quad\Longrightarrow\quad
\hhat_f\bigl(f^n(P)\bigr) > d^n\log(T) + d^n\Cr{hfh}(f)
\quad\text{for all $n\in\NN$, from \eqref{eqn:canht1},} \\
&\quad\Longrightarrow\quad
h\bigl(f^n(P)\bigr) > d^n\log(T) + (d^n-1)\Cr{hfh}(f)
\quad\text{for all $n\in\NN$, from \eqref{eqn:canht2},} \\
&\quad\Longrightarrow\quad
h\bigl(f^n(P)\bigr) > \log(T) 
\quad\text{for all $n\in\NN$.} 
\end{align*}
Hence
\begin{equation}
\label{eqn:ctfororb3}
h(P) > \log(T)+2\Cr{hfh}(f)
\quad\Longrightarrow\quad
\Count\bigl(\Orbit_f(P),T\bigr)=0.
\end{equation}
\par
The Zariski closed set~$Y$ consists of a finite number of irreducible subvarieties of~$\PP^N$ of dimension at most~$N-1$. It follows Lemma~\ref{lemma:ctYKTleCTN}(b) that
\begin{equation}
\label{eqn:countYKT}   
\Count\bigl( Y(K), T \bigr) 
\le \Cl{countPKT}(K,Y) \cdot T^N.
\end{equation}
\par
We estimate
\begin{align}
\label{eqn:contYKTorb}
\Count \biggl( & \bigcup_{P\in Y(K)} \Orbit_f(P),\, T \biggr) \notag \\
&\le \sum_{P\in Y(K)} \Count\bigl(\Orbit_f(P),T\bigr) \notag \\
&= \sum_{\substack{P\in Y(K)\\ h(P) \le \log(T)+2\Cr{hfh}(f)\\ }}
\Count\bigl(\Orbit_f(P),T\bigr) 
\quad\text{from \eqref{eqn:ctfororb3},} \notag \\
&\le \sum_{\substack{P\in Y(K)\\ h(P) \le \log(T)+2\Cr{hfh}(f)\\ }}
\Cr{Oplus}(K,N,f)\cdot\log\log(T)
\quad\text{from \eqref{eqn:ctfororb2},} \notag \\
&= \Count\bigl( Y(K), \Cl{hfh2}(f)\cdot T \bigr) \cdot \Cr{Oplus}(K,N,f)\cdot\log\log(T) \notag \\
&\omit\hfill\quad\text{where $\Cr{hfh2}=\exp\bigl(2\Cr{hfh}(f)\bigr)$,} \notag \\
&\le 
\Cr{countPKT}(K,Y) \cdot \bigl(\Cr{hfh2}(f)\cdot T\bigr)^N
\cdot \Cr{Oplus}(K,N,f)\cdot\log\log(T)
\quad\text{from \eqref{eqn:countYKT},} \notag \\
&\le \Cl{countPKTY2}(K,Y,N,f) \cdot T^N \cdot \log\log(T).
\end{align}
\par
Combining~\eqref{eqn:countPNKT} and~\eqref{eqn:contYKTorb} yields
\begin{multline*}
\Count\biggl( \PP^N(K) \setminus \bigcup_{P\in Y(K)} \Orbit_f(P),\, T \biggr) \\
\ge \Cr{PNKT}(K,N) T^{N+1} - \Cr{countPKTY2}(K,Y,N,f) \cdot T^N \cdot \log\log(T).
\end{multline*}
This function grows faster than any multiple of~$T^N$, so~\eqref{eqn:countYKT} tells us that
\[
\PP^N(K) \setminus \bigcup_{P\in Y(K)} \Orbit_f(P)
\quad\text{\parbox[t]{.4\hsize}{is not contained in a Zariski closed subset of~$\PP^N$.}}
\]
This completes the proof of~(b).
\par\noindent(c)\enspace\framebox{$\deg(f)=1$}\enspace
The assumption that~$\deg(f)=1$ tells us that~$f\in\PGL_{N+1}(\Qbar)$. Making a change of coordinates, we may assume that~$f$ is represented by a matrix~$A_f\in\GL_{N+1}(\Qbar)$ in Jordan normal form, say
\[
A_f = \begin{pmatrix}
\l_{0} & *      & 0      &\cdots & 0 \\
0      & \l_{1} & *      &\cdots & 0 \\
0      & 0      & \l_{2} &\cdots & 0 \\
\vdots & \vdots & \vdots &\ddots & \vdots \\
0       & 0     & 0      & \cdots & \l_{N} \\
\end{pmatrix} = \L + \Theta,
\]
where the stars are~$0$ or~$1$. As indicate, we let~$\L$ denote the diagonal matrix with entries~$\l_0,\ldots,\l_N$, and we let~$\Theta=A_f-\L$ be the nilpotent matrix containing the off-diagonal entries of~$A_f$. In particular, we have
\[
\L\Theta=\Theta\L\quad\text{and}\quad\Theta^N=0.
\]
It follows that for all~$n\in\ZZ$, 
\begin{equation}
\label{eqn:Afnj0N1}
A_f^n = \sum_{j=0}^{N-1} \binom{n}{j} \L^{n-j} \Theta^j.  
\end{equation}
We note that~\eqref{eqn:Afnj0N1} holds for all integers~$n$, using the usual definition of~$\binom{n}{j}$ for~$n<0$, and that it holds for~$0\le{n}<N$, since~$\binom{n}{j}=0$ for~$j>n$.
\par
We claim that our assumption that~$f$ has propagation property~(A) implies in particular that the eigenvalues are non-zero, i.e., none of the diagonal entries are~$0$. To see this, suppose the contrary. Then due to the configuration of Jordan normal form, one of the rows of~$A_f$ is~$0$, say the~$k$th row, where~$k$ is some value between~$0$ and~$N$. It follows that for any point~$P\in\PP^N$, the orbit of~$P$ is not Zariski dense; more precisely,
\[
\Orbit_f(P) \subset \{P\} \cup \{x_k=0\} \subsetneq\PP^N.
\]
This contradicts our assumption that there is at least one Zariski dense orbit, which completes the proof of the claim that the~$\l_i$ are all non-zero.
\par
We set the following notation:
\[
\begin{array}{c@{\quad}l}
K/\QQ & \text{a number field containing all of the~$\l_{i}$} \\[1\jot]
S & \text{\parbox[t]{.85\hsize}{a set of places of $K$, including all archimedean places,
such that $\l_{i}\in R_S^*$ for all $i$}} \\[5\jot]
R_S & \text{the ring of $S$-integers of $K$} \\
\end{array}
\]
For~$n\in\ZZ$ we write the corresponding power of~$A_f$ as
\[
A_f^n = \begin{pmatrix}
\mu^{(n)}_{00} & \mu^{(n)}_{01} & \cdots & \mu^{(n)}_{0N} \\
0       & \mu^{(n)}_{11} & \cdots & \mu^{(n)}_{1N} \\
\vdots  & \vdots  & \ddots & \vdots \\
0       & 0       & \cdots & \mu^{(n)}_{NN} \\
\end{pmatrix}.
\]
We note that
\[
\mu^{(n)}_{ii} = \l_{i}^n \in R_S^*~\text{for all $i$,}
\quad\text{and}\quad
\mu^{(n)}_{ij} \in R_S~\text{for all $i,j$.}
\]
\par
We are first going to prove that~$f$ has propagation property~(C2$\forall$). We let~$\Qcal\subset\PP^N(K)$ be a complete set of representatives for $\PP^N(K)/\fgeq$, and we need to show that~$\Qcal$ is Zariski dense in~$\PP^N$. So we let~$\f(\bfx)$ be a homogeneous polynomial in~$K[\bfx]$ satisfying
\begin{equation}
\label{eqn:QYinphix0}
\f(Q)=0\quad\text{for all $Q\in\Qcal$,}
\end{equation}
and our goal is to prove that~$\f=0$. 
\par
We let~$d=\deg(\f)$, and we write~$\f$ explicitly as
\[
\f(\bfx) = \sum_{k_0+k_1+\cdots+k_N=d} a(k_0,k_1,\ldots,k_N)x_0^{k_0}x_1^{k_1}\cdots x_N^{k_N}
= \sum_{|\bfk|=d} a(\bfk)\bfx^\bfk,
\]
where as indicated, we adopt the notation
\begin{equation}
\label{eqn:L1normk}
\bfk=(k_0,\ldots,k_N)
\quad\text{and}\quad
|\bfk|=k_0+k_1+\cdots+k_N.
\end{equation}
We adjoin finitely many additional primes to the set~$S$ so that the non-zero coefficients of~$\f$ are~$S$-units, i.e., 
\[
a(\bfk)\ne0 \quad\Longrightarrow\quad a(\bfk)\in R_S^*.
\]
\par
By assumption, the set of $K$-rational points~$\PP^N(K)$ is a disjoint union
\[
\PP^N(K)
= \bigsqcup_{Q\in\Qcal} \Orbit_f^\grand(Q).
\]
Hence for all $P\in\PP^N(K)$ there exists a point $Q_P\in\PP^N(K)$ and an integer $n_P\in\ZZ$ such that\footnote{Since~$n_P\in\ZZ$, it's notationally convenient to insert a negative sign here.}
\[
\f(Q_P)=0 \quad\text{and}\quad P = f^{-n_P}(Q_P).
\]
Note that we have complete freedom in our choice of~$P\in\PP^N(K)$. 
\par
We fix a prime~$\gp$ of~$R_K$ with~$\gp\notin{S}$, and we let~$\pi\in{R_K}$ be a uniformizer for~$\gp$.  Then for each 
\[
\bfr=(r_0,r_1,\ldots,r_N)\in\NN^{N+1},
\quad\text{we set}\quad
P(\bfr) = [\pi^{-r_0},\pi^{-r_1},\ldots,\pi^{-r_N}].
\]
For the moment we  assume that the~$r_i$ are positive and decreasing,
\begin{equation}
\label{eqn:r0gtr1gtrNge1}
r_0 > r_1 > \cdots > r_N \ge 1.
\end{equation}
Later we will specify them more precisely.
\par
We let~$0\le{i}\le{N}$, and for~$n\in\ZZ$ we consider the~$\gp$-adic valuation of\footnote{More formally, we lift~$P(\bfr)$ to~$(\pi^{-r_0},\ldots,\pi^{-r_N})\in\AA^{N+1}$, apply~$A_f^n$, and take the~$i$th coordinate; but for ease of exposition, we will simply refer to the~$i$th coordinate in~$\PP^N$.}
\[
f^n\bigl(P(\bfr)\bigr)[i]
:= \text{the $i$th coordinate of~$f^n\bigl(P(\bfr)\bigr)$.}
\]
Thus
\begin{align*}
\ord_\gp\Bigl(f^n\bigl(P(\bfr)\bigr)[i]\Bigr)
&=
\ord_\gp\biggl( \sum_{j=i}^N \mu_{ij}^{(n)} \cdot \pi^{-r_j} \biggr) \\
&=
\ord_\gp\biggl( \l_i^n \pi^{-r_i} + \sum_{j=i+1}^N \mu_{ij}^{(n)} \cdot \pi^{-r_j} \biggr).
\end{align*}
The facts that~$\gp\notin{S}$ and~$\l_i\in{R_S^*}$ and $\mu_{ij}^{(n)}\in{R_S}$ yield
\[
\ord_\gp( \l_i^n \pi^{-r_i} ) = -r_i
\quad\text{and}\quad
\ord_\gp( \mu_{ij}^{(n)} \cdot \pi^{-r_j} ) \ge -r_j.
\]
It follows from~\eqref{eqn:r0gtr1gtrNge1} that we have a strict inequality
\[
\ord_\gp( \l_i^n \pi^{-r_i} ) < \ord_\gp( \mu_{ij}^{(n)} \cdot \pi^{-r_j} )
\quad\text{for all $i<j$,} 
\]
and hence the non-archimedean triangle inequality yields
\begin{equation}
\label{eqn:ordpfnPri}
\ord_\gp\Bigl(f^n\bigl(P(\bfr)\bigr)[i]\Bigr) 
= \ord_\gp( \l_i^n \pi^{-r_i} ) = -r_i.
\end{equation}
\par
We compute
\begin{align}
\label{eqn:0fQPr}
0 = \f(Q_{P(\bfr)}) 
&= \f\Bigl(f^{n_{P(r)}} \bigl(P(\bfr)\bigr) \Bigr)  \notag\\
&= \sum_{|\bfk|=d} a(\bfk) \Bigl( f^{n_{P(\bfr)}}\bigl(P(\bfr)\bigr) \Bigr)^\bfk \notag\\
&= \sum_{|\bfk|=d} a(\bfk) \prod_{i=0}^N f^{n_{P(\bfr)}}\bigl(P(\bfr)\bigr)[i]^{k_i} .
\end{align}
\par
We use~\eqref{eqn:ordpfnPri} to compute the valuation of the non-zero monomials appearing in~\eqref{eqn:0fQPr}. Thus if~$a(\bfk)\ne0$, then
\begin{align}
\label{eqn:ordpaiLnPpi2}
\ord_\gp\Bigl( a(\bfk) &  \Bigl( f^{n_{P(\bfr)}}\bigl(P(\bfr)\bigr) \Bigr)^\bfk \notag\\*
&= 
\ord_\gp\bigl(a(\bfk)\bigr) + 
\sum_{i=0}^N  \ord_\gp\Bigl( f^{n_{P(\bfr)}}\bigl(P(\bfr)\bigr)[i]^{k_i} \Bigr) \notag\\*
&=
\sum_{i=0}^N -r_i k_i.
\end{align}
The last equality follows from~\eqref{eqn:ordpfnPri} and the fact that the non-zero~$a(\bfk)$ are~$S$-units.
\par
We now set the~$r_i$ to equal
\[
r_i = (d+1)^{N-i} \quad\text{for $0\le i\le N$.}
\]
It then follows from~\eqref{eqn:ordpaiLnPpi2} and Lemma~\ref{lemma:kmapstordotkinjective} that the non-zero monomials appearing in the expansion~\eqref{eqn:0fQPr} have distinct negative $\gp$-adic valuations. Thus if~$\f(\bfx)$ has any non-zero monomials, then the non-archimedean triangle inequality implies that the sum cannot equal~$0$. Since $\f(Q_{P(\bfr)})=0$ by construction, it follows that the polynomial~$\f$ has no non-zero monomials, i.e., we have proven that~$\f=0$, which completes the proof that~$\Qcal$ is Zariski dense in~$\PP^N$, and thus the proof that the map~$f$ has property~(C2$\forall$).
\par
We next invoke Lemma~\ref{lemma:linmapXfdenseeqtorus}, which tells us that the set of dense-orbit points~$(\PP^N)_f^\dense$ contains a non-empty Zariski open subset of~$\PP^N$. This fact, combined with the assumed propagation property~(A) and the proven propagation property~(C2$\forall$) allows us to apply Lemma~\ref{lemma:C2forallZDimpliesC3forall} and conclude that~$f$ has propagation property~(C3$\forall$).
\end{proof}

The following result, which is used in the proof of Theorem~\ref{theorem:projectivespacedeg2}(c), is essentially proven by Ghioca and Hu in~\cite{MR3829742}. We indicate how to modify their proof to obtain the desired result.

\begin{lemma}
\label{lemma:linmapXfdenseeqtorus}
Let $f:\PP^N\to\PP^N$ be a linear map, and assume that there is at least one Zariski-dense~$f$-orbit. Then~$X_f^\dense$ contains a non-empty Zariski open subset of~$\PP^N$. 
\end{lemma}
\begin{proof}
By assumption, there is a point~$P_0\in\PP^N$ whose~$f$-orbit is Zariski dense. We change coordinate so that~$f$ is represented by a matrix in Jordan normal form. Then the hyperplane~$H:=\{X_N=0\}$ is~$f$-invariant, so letting \text{$\AA^N:=\PP^N\setminus H_N$} and dehomogenizing, we see that~$f$ induces a linear map on
\[
A : \AA^N\longrightarrow \AA^N
\quad\text{given by a matrix}\quad
A\in\GL_N(K).
\]
Further, we must have~$P_0\in\AA^N(K)$, since~$f(H)\subseteq{H}$, so points in~$H_N$ do not have Zariski dense orbits.
\par
We now follow the proof of~\cite[Theorem~2.1]{MR3829742}. The assumption that there is a Zariski-dense orbit implies in particular that there are no non-constant rational functions~$\f:\AA^N\dashrightarrow\PP^1$ satisfying~\text{$\f\circ{A}=\f$}. (Otherwise $A$-orbits would lie in the fibers of~$\f$.) The conclusion of~\cite[Theorem~2.1]{MR3829742} is that there exists at least one Zariski-dense $A$-orbit in~$AA^N$, but the proof actually shows that the~$A$-orbit of the specific point~$(1,1,\ldots,1)$ is Zariski dense in~$\AA^N$. (Remember that we are assuming that~$A$ is in Jordan normal form.) We briefly indicate how to modify the proof of~\cite[Theorem~2.1]{MR3829742} to show that every point
\[
\bfbeta=(\b_1,\ldots,\b_N)
\in \TT^N 
:= \bigl\{[x_0,x_1,\ldots,x_N]\in\PP^N:x_0\cdots x_N\ne0\bigr\}
\]
has Zariski-dense~$A$-orbit.
\par
The proof in~\cite[Theorem~2.1]{MR3829742} that~$(1,1,\ldots,1)$ has Zariski-dense~$A$-orbit is reduced to two cases, labeled Case~4 and~5 in~\cite{MR3829742}. We start with Case~4, which is the case that~$A$ is a diagonal matrix with multiplicatively independent diagonal entries.\footnote{We remark that in Case~4 there is an easier proof that avoids the use of deep results on linear recurrences and shows that $(\PP^N)_f^\dense=\TT^N$. One simply starts with one point with dense orbit and translates it to all of~$\TT^N$ using diagonal matrices, while noting that diagonal matrices commute with~$A$ and thus preserve the Zariski density property of the orbit. However, it is unclear whether this translation proof can be adapted to handle the non-diagonal matrices in Case~5.} Let
\begin{equation}
\label{eqn:Fx1xNci1iN}
F(x_1,\ldots,x_N) = \sum_{i_1,\ldots,i_N} c_{i_1,\ldots,i_N}\prod_{j=1}^N x_j^{i_j}
\end{equation}
be a polynomial such that~$F(A^n\bfbeta)=0$ for all~$n\ge0$. Our goal is to show that~$F=0$. Letting \text{$\L_{i_1,\ldots,i_N}=\l_1^{i_1}\cdots\l_N^{i_N}$}, we have
\begin{equation}
\label{eqn:FAnbfbetaeq0}
F(A^n\bfbeta)
= \sum_{i_1,\ldots,i_N}
\left(c_{i_1,\ldots,i_N}\b_1^{i_1}\cdots\b_N^{i_N}\right)\cdot\L_{i_1,\ldots,i_N}.
\end{equation}
As explained in~\cite{MR3829742}, the sequence~\eqref{eqn:FAnbfbetaeq0} is a non-degenerate linear recurrence as a function of~$n$, so it can vanish for only finitely many~$n$ unless all of the  coefficients~$c_{i_1,\ldots,i_N}\b_1^{i_1}\cdots\b_N^{i_N}$ vanish. But the~$\b_j$ are non-zero by assumption, so this forces the~$c_{i_1,\ldots,i_N}$ to vanish, which completes the proof that~$F=0$.
\par
Similarly, in Case~5 of~\cite[Theorem~2.1]{MR3829742} the matrix~$A$ is diagonal except for one $2$-dimensional Jordan block. So letting~$\l_1=\l_2$ be the eigenvalue in the Jordan block, we have
\[
A^n(\bfbeta) =(\b_1\l_1^n+\b_2n\l_1^{n-1},\b_2\l_1^n,\b_3\l_3^n,\ldots,\b_N\l_N^n).
\]
We define a linear transformation
\[
B(x_1,\ldots,x_N) = \bigl( \l_1(\b_2x_1-\b_1x_2),x_2,\ldots,x_N\bigr).
\]
The fact that~$\b_2\ne0$ ensures that~$B$ is invertible, so it suffices to show that~$B\Orbit_A(\bfbeta)$ is Zariski dense. We have
\[
BA^n(\bfbeta) = (\b_2^2n\l_1^n,\b_2\l_1^n,\b_3\l_3^n,\ldots,\b_N\l_N^n).
\]
Let~$F(\bfx)$ be a polynomial~\eqref{eqn:Fx1xNci1iN} vanishing on~$B\Orbit_A(\bfbeta)$. Then
\begin{align*}
F\bigl(BA^n(\bfbeta)\bigr) 
&= \sum_{i_1,\ldots,i_N} c_{i_1,\ldots,i_N}
(\b_2^2n\l_1^n)^{i_1} \cdot (\b_2\l_1^n)^{i_2}
\cdot \prod_{j=3}^N (\b_j\l_j^n)^{i_j} \\
&= \sum_{k,i_3,\ldots,i_N}
\biggl(\sum_{i_1+i_2=k }
c_{i_1,\ldots,i_N} \b_2^{2i_1+i_2}\b_3^{i_3}\cdots\b_N^{i_N}n^{i_1}\biggr)\l_1^k\l_3^{i_3}\cdots\l_N^{i_N}.
\end{align*}
As in the proof of Case~5 of~\cite[Theorem~2.1]{MR3829742}, the sequence~$F\bigl(BA^n(\bfbeta)\bigr)$ is a non-degenerate linear recurrence, so if it vanishes for infinitely many~$n\ge0$ (much less for all $n\ge0$), every coefficient must vanish, i.e., 
\[
c_{i_1,\ldots,i_N} \b_2^{2i_1+i_2}\b_3^{i_3}\cdots\b_N^{i_N} = 0
\quad\text{for all $i_1,\ldots,i_N$.}
\]
Then using the fact that~$\b_2\cdots\b_N\ne0$, we conclude that~$F=0$.
\par
This concludes the proof that after an appropriate change of coordinates, the torus~$\TT^N\subset\AA^N\subset\PP^N$ is contained in the set~$(\PP^N)_f^\dense$ of points with Zariski-dense~$f$-orbits. 
\end{proof}

We conclude this section with two elementary results.  The first describes intervals contained within other intervals that was used in the proof of Theorem~\ref{theorem:projectivespacedeg2}(a). The second is an injectivity result for dot products of sequences.

\begin{lemma}
\label{lemma:missedinterval}
Let $\a_1,\ldots,\a_r>0$ and $\b\ge0$ and $d>1$. For~$1\le{i}\le{r}$ and~$n\in\NN$, define real intervals
\[
I_{i,n} := [\a_i d^n - \b\,, \a_i d^n + \b]
\quad\text{and}\quad I := \bigcup_{1\le i\le r} \bigcup_{n\in\NN} I_{i,n}.
\]
There are constants $\Cl{aibde}$ and $\Cl{aibde7}>0$, depending only on~$\a_1,\ldots,\a_r,\b,d$, so that
\[
T \ge \Cr{aibde}
\quad\Longrightarrow\quad
\left(\begin{array}{@{}l@{}}
\text{$[0,T]\setminus I$ contains an interval} \\
\text{of length at least ${\Cr{aibde7}T}/{\log(T)}$.}
\end{array}\right).
\]
\end{lemma}
\begin{proof}
Fix~$T$. In the following calculation, all constants are independent of~$T$. The number of intervals~$I_{i,n}$ that have a point in common with the interval~$[0,T]$ is bounded by
\begin{align*}
\#\bigl\{(i,n)\in[r]\times\NN & : I_{i,n}\cap[0,T]\ne\emptyset \bigr\}\\
&\le \sum_{i=1}^r \#\bigl\{n \in \NN : \a_id^n-\b \le T \bigr\} \\
&\le \sum_{i=1}^r \left( \log_d\left(\frac{T+\b}{\a_i}\right) + 2 \right) \\
&\le \Cl{aibde2}\log_d(T) \quad\text{for $T\ge\Cl{aibde3}$.}
\end{align*}
Suppose that $[0,T]\setminus{I}$ contains no intervals of length at least~$B$. This implies that if we lengthen each~$I_{i,n}$ by~$B/2$ on each side, then~$\{I_{i,n}\}$ will cover all of~$[0,T]$. In other words, if we define
\[
I_{i,n}(B) := [\a_i d^n - \b - B/2\,, \a_i d^n + \b + B/2],
\]
then the assumption that~$[0,T]\setminus{I}$ contains no intervals of length at least~$B$ implies that
\[
[0,T] \subseteq \bigcup_{\substack{i\in[r],\,n\in\NN\\ I_{i,n}\cap[0,T]\ne\emptyset\\}} I_{i,n}(B). 
\]
Hence
\begin{align*}
T 
&\le \sum_{\substack{i\in[r],\,n\in\NN\\ I_{i,n}\cap[0,T]\ne\emptyset\\}} 
\operatorname{Length}\bigl(I_{i,n}(B)\bigr) \\
&\le \sum_{\substack{i\in[r],\,n\in\NN\\ I_{i,n}\cap[0,T]\ne\emptyset\\}} 
(2\b+B) \\
&\le (2\b+B)\cdot \bigl( \Cr{aibde2} \log_d(T) \bigr)
\quad\text{for $T\ge\Cr{aibde3}$.} 
\end{align*}
Rearranging this inequality, we have proven that if~$T\ge\Cr{aibde3}$ and if $[0,T]\setminus{I}$ contains no intervals of length at least~$B$, then
\[
B \ge \frac{T}{\Cr{aibde2}\log_d(T)} - 2\b.
\]
Adjusting constants, it follows that if~$T\ge\Cl{aibde4}$, then~$[0,T]\setminus{I}$ will contain an interval of length $\Cl{aibde5}{T}/\log(T)$.
\end{proof}

\begin{lemma}
\label{lemma:kmapstordotkinjective}
Define a list of integers
\begin{equation}
\label{eqn:rjdj1d1x}
r_i = (d+1)^{N-i} \quad\text{for $0\le i\le N$.} 
\end{equation}
Then with notation as in~\eqref{eqn:L1normk}, the map
\begin{equation}
\label{eqn:idtorix}
\bigl\{ \bfk\in\NN^{N+1} : |\bfk|=d \bigr\} \longrightarrow \NN,\quad
\bfk \longmapsto \bfr\cdot\bfk =\sum_{i=0}^N r_i k_i
\end{equation}
is injective.
\end{lemma}
\begin{proof}
We suppose that 
\[
\bfj\ne\bfk 
\quad\text{and}\quad
|\bfj|=|\bfk|=d
\quad\text{and}\quad
\bfr\cdot\bfj=\bfr\cdot\bfk
\]
and derive a contradiction. The assumption that~$\bfj\ne\bfk$ implies that there is an index~$0\le\ell\le{N}$ such that
\begin{equation}
\label{eqn:iellnekell2}
j_\ell\ne k_\ell \quad\text{and}\quad j_m=k_m~\text{for $0\le m<\ell$.}
\end{equation}
We now compute 
\begin{align*}
\bfr\cdot\bfj &=\bfr\cdot\bfk \\
&\quad\Longrightarrow\quad
\sum_{i=\ell}^N r_i j_i = \sum_{i=\ell}^N r_i k_i
\quad\text{from \eqref{eqn:iellnekell2},} \\
&\quad\Longrightarrow\quad
\bigl| r_\ell (j_\ell-k_\ell) \bigr|
\le \sum_{i=\ell+1}^N \bigl| r_i (j_i-k_i) \bigr|
\quad\text{triangle inequality,}\\
&\quad\Longrightarrow\quad
(d+1)^{N-\ell}  |j_\ell-k_\ell|
\le \sum_{i=\ell+1}^N (d+1)^{N-i} |j_i-k_i| 
\quad\text{from \eqref{eqn:rjdj1d1x},} \\
&\quad\Longrightarrow\quad
(d+1)^{N-\ell} \le \sum_{i=\ell+1}^N (d+1)^{N-i} |j_i-k_i|  
\quad\text{since $j_\ell\ne k_\ell$ from \eqref{eqn:iellnekell2},} \\
&\quad\Longrightarrow\quad
(d+1)^{N-\ell} \le (d+1)^{N-\ell-1} \cdot d + (d+1)^{N-\ell-2} \cdot d \\
&\omit\hfill\text{\parbox{.6\hsize}{since $|\bfj|=|\bfk|=d$, so the right-hand side is 
maximized by taking $j_{\ell+1}=k_{\ell+2}=d$,}} \\
&\quad\Longrightarrow\quad
(d+1)^2 \le (d+1)\cdot d + d = (d+1)^2-1 \\
&\omit\hfill\text{\parbox{.6\hsize}{dividing both sides by $(d+1)^{N-\ell-2}$.}}
\end{align*}
This contradiction proves that the map~\eqref{eqn:idtorix} is injective for the choice~\eqref{eqn:rjdj1d1x} of~$r_0,\ldots,r_N$. 
\end{proof}

\section{K3 Surfaces}
\label{section:K3surfaces}

\begin{theorem}
\label{theorem:K3}
Let~$X/\Qbar$ be a smooth projective~K3 surface, let~$f:X\to{X}$ be an endomorphism, and assume that there is a point~$P_0\in{X(\Qbar)}$ whose $f$-orbit~$\Orbit_f(P_0)$ is Zariski dense in~$X$. Then there is a common field of definition~$K$ for~$X$,~$f$ and~$P_0$ such that for any finite collection of grand $f$-orbits $\G_1,\ldots,\G_r\subset{X}(\Qbar)$, we have
\begin{equation}
\label{eqn:K3XminusGi}
\overline{X(K) \setminus (\G_1\cup\cdots\cup\G_r)} = X.
\end{equation}
In the terminology of Table~\textup{\ref{table:propogationstatements}}, endomorphisms of~$X$ satisfy the orbit propagation statement
\[
\textup{(A)} \quad\Longrightarrow\quad \textup{(C1)}.
\]
\end{theorem}
\begin{proof}
We first recall the well-known fact that all dominant endomorphisms of a smooth algebraic~K3 surface~$X$ are automorphisms; cf.~\cite[Section1.2]{MR2597304}. To see this, we start with the general formula
\[
f^*K_X = K_X + \text{(Ramification divisor of $f$)}.
\]
Since~K3 surfaces have trivial canonical bundle~$K_X=0$, we conclude that~$f$ is unramified. But~$X(\CC)$ is simply connected, and hence~$f$ is an automorphism.
\par
We note that the existence of a point with Zariski dense~$f$-orbit implies in particular that~$f$, and hence also~$f^{-1}$, have infinite order.
Further, as noted in Lemma~\ref{lemma:orbitelemproperties}, the fact that~$f$ is an automorphism implies that the grand $f$-orbit of a point~$P$ is the union of the $f$-orbit of~$P$ and the~$f^{-1}$-orbit of~$P$. Hence it suffices to prove~\eqref{eqn:K3XminusGi} for~$f$ under the assumption that the~$\G_i$ are forward $f$-orbits, since we can then apply the same reasoning to~$f^{-1}$. 
\par
\par
We say that a curve~$C\subset{X}$ is \emph{$f$-periodic} if~$f^n(C)={C}$ for some $n\ge1$.\footnote{We note that since~$f$ is an automorphism, a curve is $f$-periodic if and only if it is~$f$-preperiodic.} Otherwise we say that~$C$ is
\emph{$f$-wandering}. Our first step is to produce an $f$-wandering rational curve on~$X$. To do this, we use the following two results.
\begin{lemma}
\label{lemma:21}
\cite[Theorem~2.1]{MR2917181}
Let~$D$ be a non-trivial effective divisor on~$X$. Then~$D$ is linearly equivalent to a sum of rational curves.
\end{lemma}
\begin{lemma}
\label{lemma:22}
\cite[Chapter~5,~Corollary~3.5]{MR3586372}
Let~$D$ be an ample divisor on~$X$, and let~$g:X\to{X}$ be an automorphism such that $g_*D\sim{D}$. Then~$g$ is of finite order.
\end{lemma}
Resuming the proof of Theorem~\ref{theorem:K3}, we suppose that all rational curves in~$X$ are $f$-periodic and derive a contradiction. Let~$D$ be an ample effective divisor on~$X$. By Lemma~\ref{lemma:21}, there are rational curves~$Y_1,\ldots,Y_m$ and integers~$a_1,\ldots,a_m$ such that
\[
D \sim a_1Y_1+\cdots+a_mY_m.
\]
The rational curves~$Y_j$ are assumed to be $f$-periodic. Let $p\ge1$ be a positive integer such that $f^p(Y_j)=Y_j$ for all~$1\le{j}\le{m}$. Since~$f$ is an automorphism, this implies that $f_*^pY_j=Y_j$ as divisors. Hence
\[
f_*^pD \sim a_1f_*^pY_1+\cdots+a_mf_*^pY_m=a_1Y_1+\cdots+a_mY_m\sim D.
\]
Since~$D$ is ample by assumption, Lemma~\ref{lemma:22} implies that~$f^p$ has finite order, which contradicts the assumption that there exists a Zariski dense $f$-orbit.
\par
We let~$Y$ be an $f$-wandering rational curve on~$X$. Replacing~$K$ with a finite extension, which by abuse of notation we continue to call~$K$, we may assume that~$Y$ is defined over~$K$ and that~$Y(K)$ is Zariski dense in~$Y$.
\par
Let~$\Fcal\subset{X(K)}$ be a finite set of points, and let
\[
\Orbit_f^\grand(\Fcal) := \bigcup_{Q\in\Fcal} \Orbit_f^\grand(Q)
= \bigcup_{Q\in\Fcal} \Bigl( \Orbit_f(Q) \cup \Orbit_{f^{-1}}(Q) \Bigr)
\]
be the union of the $f$-orbits of the points in~$\Fcal$. 
Our goal is to prove that~$X(K)$ has a Zariski dense set that is disjoint from~$\Orbit_f^\grand(\Fcal)$. 
\par
We know from earlier that~$Y(K)$ is Zariski dense in~$Y$. Further, the dynamical Mordell--Lang conjecture is true for \'etale maps~\cite{MR2766180}, so in particular for automorphisms such as~$f$ and~$f^{-1}$.  Hence 
\begin{equation}
\label{eqn:MLforetale}
\left(\begin{tabular}{@{}l@{}}
$\Orbit_f^\grand(Q)\cap C$ is finite for every point \\
$Q\in X(K)$ and every  curve~$C\subset{X}$ \\
\end{tabular}
\right).
\end{equation}
In particular, since~$f^n(Y)$ is a curve for every~$n\ge0$, and since~$\Orbit_f^\grand(\Fcal)$ is a finite union of $f$-orbits and~$f^{-1}$-orbits, it follows
from~\eqref{eqn:MLforetale} that
\begin{equation}
\label{eqn:fnRZDinX2}
\underbrace{f^n\bigl(Y(K)\bigr)\setminus\Orbit_f^\grand(\Fcal)}_{\hidewidth\text{This is an infinite subset of the curve $f^n(Y)$.}\hidewidth}\quad\text{is Zariski dense in $f^n(Y)$.}
\end{equation}
But the fact that~$Y$ is wandering implies that the union of curves
\begin{equation}
\label{eqn:fnRZDinX3}
\bigcup_{n\ge0} f^n(Y)\quad\text{is Zariski dense in $X$,}
\end{equation}
and combining~\eqref{eqn:fnRZDinX2} and~\eqref{eqn:fnRZDinX3}, we conclude that 
\[
\bigcup_{n\ge0} \Bigl( f^n\bigl(Y(K)\bigr)\setminus\Orbit_f^\grand(\Fcal) \Bigr)
\quad\text{is Zariski dense in $X$.}
\]
Hence
\[
X(K) \setminus \Orbit_f^\grand(\Fcal)
\supseteq 
\biggl( \bigcup_{n\ge0} \Bigl( f^n\bigl(Y(K)\bigr) \biggr)
 \setminus \Orbit_f^\grand(\Fcal),
\]
is also Zariski dense in~$X$, which completes the proof that~(A) implies~(C1) for K3 surfaces.
\end{proof}

\section{Abelian Varieties}
\label{section:abelianvariety}

In this section we prove that our strongest propagation principle is true for geometrically simple abelian varieties.

\begin{theorem}
\label{theorem:simpleabelianvariety}
Let~$X/\Qbar$ be a geometrically simple abelian variety, let~$f:X\to{X}$ be an endomorphism of~$X$ \textup(as an abstract variety\textup), and let~$P_0\in{X}(\Qbar)$ be a point such that~$\Orbit_f(P_0)$ is Zariski dense in~$X$. Then there is a number field~$K/\QQ$ such that~$X$ and~$f$ are defined over~$K$ and such that every complete set of representatives for 
\[
\text{$X_f^\dense(K)/\fgeq$\quad  is Zariski dense in $X$.}
\]
In the terminology of Table~\textup{\ref{table:propogationstatements}}, endomorphisms of geometrically simple abelian varieties satisfy the orbit propagation statement
\[
\textup{(A)} \quad\Longrightarrow\quad \textup{(C3$\forall$)}.
\]
\end{theorem}

\begin{remark}
Let~$X$ be a split semi-abelian variety, and let $f:X\to{X}$ an endomorphism, all defined over~$\Qbar$. There is an explicit construction in~\cite{MR4396628} that produces points in~$X(\Qbar)$ having Zariski dense $f$-orbits. The paper does not appear to discuss whether one can create many such points defined over a fixed number field, but the construction is sufficiently explicit that it might give an alternative approach to proving orbit propagation results for (semi)-abelian varieties.
\end{remark}

\begin{proof}[Proof of Theorem \ref{theorem:simpleabelianvariety}]
We start by fixing a field of definition~$K$ for~$X$,~$f$, and~$P_0$, but we may at times replace~$K$ with a finite extension. Every map between abelian varieties is the composition of a homomorphism and a translation~\cite[Section~4, Corollary~1]{MR2514037}, say
\[
f(x) = \f(x) + Q_0,\;
\text{with $Q_0\in X(K)$ and $\f:X\to X$ a homomorphism.}
\]
The assumption that~$P_0\in{X}(K)$ satisfies
\[
\overline{\Orbit_f(P_0)} = X
\]
implies, in particular, that~$f$ is dominant, and thus that~$\f$ is a non-zero isogeny. 
\par
The assumption that~$X$ is geometrically simple implies that the  kernel of~$\f-1$ is either a finite subgroup of~$X$ or all of~$X$. We consider these two case in turn.
\Case{1}{$\Ker(\f-1)$ is finite}
In this case~$\f-1$ is an isogeny, i.e., the map~$\f-1:X\to{X}$ is a finite surjective map. Hence we can find a point~$R_0\in{X(\Kbar)}$ satisfying
\[
(\f-1)(R_0) = Q_0.
\]
Writing~$T_P$ in general for the translation-by-$P$ map, we have
\begin{align*}
T_{R_0}\circ f\circ T_{-R_0}(x)
&= f(x-R_0) + R_0\\
&= \f(x-R_0)+Q_0+R_0\\
&= \f(x)-\f(R_0)+Q_0+R_0\\
&= \f(x).
\end{align*}
Thus~$f$ and~$\f$ are conjugate via the automorphism~$T_{R_0}\in\Aut(X)$, so they have the same dynamics. Hence in Case~1, it suffices to consider the case that~$f$ itself is an isogeny.
\par
Let~$\Mcal$ be an infinite set of positive integers with the following properties:
\begin{align}
\label{eqn:Mgcdfm1}
\gcd\bigl(\deg(f),m\bigr)&=1~\text{for all $m\in\Mcal$.} \\
\label{eqn:Mgcdm1m2}\gcd(m_1,m_2)&=1~\text{for all distinct $m_1,m_2\in\Mcal$.} 
\end{align}
For example, we could take~$\Mcal$ to be the set of all primes not dividing~$\deg(f)$. Then for~$m_1,m_2\in\Mcal$ we compute
\begin{align*}
    \Orbit_f&(m_1P_0)  \cap \Orbit_f(m_2P_0) \ne \emptyset \\
    &\quad\Longleftrightarrow\quad
    f^{n_1}(m_1P_0) = f^{n_2}(m_2P_0) \quad\text{for some $n_1,n_2\in\NN$,} \\
    &\quad\Longrightarrow\quad
    (f^{n_1}\circ m_1 -  f^{n_2}\circ m_2)(P_0)=0 \quad\text{for some $n_1,n_2\in\NN$,} \\
    &\quad\Longrightarrow\quad
   f^{n_1}\circ m_1 -  f^{n_2}\circ m_2 = 0 \quad
   \text{\parbox[t]{.4\hsize}{in $\ZZ[f]\subseteq\Isog(X)$ for some $n_1,n_2\in\NN$,
   since $P_0$ is non-torsion because its $f$-orbit is Zariski dense,}} \\
    &\quad\Longrightarrow\quad
   (\deg f)^{n_1}\cdot m_1^{2g} = \deg(f)^{n_2}\cdot m_2^{2g} 
   \quad\text{for some $n_1,n_2\in\NN$,} \\
    &\quad\Longrightarrow\quad
    m_1=m_2 \quad\text{from \eqref{eqn:Mgcdfm1} and \eqref{eqn:Mgcdm1m2}.}
\end{align*}
Taking the contrapositive, we have proven that
\[
\text{$m_1,m_2\in\Mcal_f$ and $m_1\ne m_2$}
\quad\Longrightarrow\quad
\Orbit_f(m_1P_0)  \cap\Orbit_f(m_2P_0)=\emptyset.
\]
Hence the grand $f$-orbits generated by the points in~$\{mP_0:m\in\Mcal_f\}$ are distinct.
\par
We claim that the $f$-orbits of these points are also Zariski dense. To prove this, we note that since~$f$ is an isogeny, it commutes with multiplication-by-$m$, and since~$f$ is a finite map of a proper variety, it sends closed sets to closed sets. Hence
\begin{align}
\label{eqn:OfmPZC}
\overline{\Orbit_f(mP_0)} &= \overline{m\cdot\Orbit_f(P_0)}
\quad\text{since $fm=mf$,} \notag\\
&= m\cdot\overline{\Orbit_f(P_0)}
\quad\text{since $m$ is finite and $X$ proper.}
\end{align}
It follows immediately from~\eqref{eqn:OfmPZC} that\footnote{The converse also holds, since if~$\Orbit_f(P_0)$ is not Zariski dense, then its Zariski closure is a subvariety whose dimension is strictly smaller than~$\dim(X)$, and multiplication-by-$m$ preserves the dimension, so~$\Orbit_f(mP_0)$ is also not Zariski dense.}
\[
\overline{\Orbit_f(P_0)} = X\quad\Longrightarrow\quad \overline{m\Orbit_f(P_0)} = X.
\]
Since the $f$-orbit of~$P_0$ is Zariski dense by assumption, the same is true of the~$f$-orbit of~$mP_0$ for all positive integers~$m$.
\par
We let
\[
\Mcal\cdot P_0 := \{mP_0:m\in\Mcal\},
\]
and we consider the following statements:
\begin{parts}
\Part{(1)}
Every point in $\Mcal\cdot{P_0}$ has a Zariski dense $f$-orbit.
\Part{(2)}
The points in $\Mcal\cdot{P_0}$ generate distinct grand $f$-orbits.
\Part{(3)}
The set $\Mcal\cdot{P_0}$ is Zariski dense.
\end{parts}
We have proven~(1) and~(2), and we now prove~(3). We first show that~$\Mcal\cdot{P_0}$ is an infinite set. Since~$\Mcal$ is an infinite set of integers, it suffices to show that~$P_0$ is a non-torsion point. To see this, we assume that~$mP_0=0$ for some non-zero integer and derive a contradiction. The fact that~$f$ is an isogeny implies that~$f$ commutes with the multiplication-by-$m$ map, so we have
\[
\Orbit_f(P_0) \subseteq \ZZ[f]\cdot P_0 = (\ZZ[f]/m\ZZ[f])\cdot P_0.
\]
The ring~$\ZZ[f]/m\ZZ[f]$ is finite, so the $f$-orbit~$\Orbit_f(P_0)$ is finite, i.e., the point~$P_0$ is~$f$-preperiodic. This contradicts the assumption that~$\Orbit_f(P_0)$ is Zariski density, and completes the proof that~$\Mcal\cdot{P_0}$ is an infinite set of points.
\par
The Zariski closure $\overline{\Mcal\cdot{P_0}}$ contains the infinite set~$\Mcal\cdot{P_0}$, which in turn is a subset of the finitely generated (indeed, rank~$1$) subgroup~$\ZZ\cdot{P_0}$ of~$X$. It follows from Faltings's theorem~\cite{MR1109353} (originally the Mordell--Lang conjecture) that~$\overline{\Mcal\cdot{P_0}}$ contains a translate of an abelian subvariety of~$X$, necessarily positive dimensional since~$\#(\Mcal\cdot{P_0})=\infty$. The assumed simplicity of~$X$ tells us that the only such abelian subvariety is~$X$ itself. Hence $\overline{\Mcal\cdot{P_0}}=X$, which completes the proof of~(3).
\par
We can restate~(1) as the inclusion~$\Mcal\cdot{P_0}\subseteq{X_f^\dense}$, and then~(3) tells us that~$X_f^\dense(K)$ contains a Zariski dense set of points. We have thus proven the following useful fact:
\begin{equation}
\label{eqn:Xfdenseisdenseabvar}
\overline{X_f^\dense(K)} = X.
\end{equation}
\par
Our next goal is to prove that~$f$ has propagation property~(C2$\forall$). So we let $\Qcal\subset{X(K)}$ be a set of representatives for the grand $f$-orbits, i.e., a set satisfying
\begin{equation}
\label{eqn:QbijectXKfgeqabvar}
\Qcal \xleftrightarrow{\;\text{bijective}\;} X(K)/\fgeq,
\end{equation}
and we need to show that~$\Qcal$ is Zariski dense. We let
\[
Y = \overline\Qcal \subseteq X,
\]
and our goal is to prove that~$Y=X$.
\par
We proved earlier that the set~$\Mcal\cdot{P_0}$ is Zariski dense and that each point in~$\Mcal\cdot{P_0}$ is in a distinct grand~$f$-orbit. The choice~\eqref{eqn:QbijectXKfgeqabvar} of~$\Qcal$ says that every grand $f$-orbit containing a~$K$-rational point will contain a unique point of~$\Qcal$. This allows us to  define an injection
\[
\psi : \Mcal\cdot{P_0} \longhookrightarrow\Qcal,
\quad
\left(
\text{\begin{tabular}{@{}l@{}}
$\psi(mP_0)$ is the unique point~$Q\in\Qcal$\\
such that $\Orbit_f(Q)^\grand=\Orbit_f(mP_0)^\grand$\\
\end{tabular}}
\right).
\]
The existence of Zariski dense $f$-orbit~$\Orbit_f(P_0)$ tells us, in particular, that~$f\ne0$, so there exists an isogeny $g:X\to{X}$ satisfying
\[
f\circ g = g\circ f = q \in \ZZ_{>0}.
\]
The condition that $\Orbit_f(Q)^\grand=\Orbit_f(mP_0)^\grand$ with $Q=\psi(mP_0)$
says that there are integers~$n_1,n_2\ge0$ satisfying 
\begin{equation}
\label{fn1Qfn2mP0abvar}
f^{n_1}(Q)=f^{n_2}(mP_0)
\end{equation}
Applying~$g^{n_1}$ to both sides of~\eqref{fn1Qfn2mP0abvar}, we find that for $mP_0\in\Mcal$ there are non-negative integers~$n_1$ and~$n_2$ satisfying
\[
q \psi(mP_0) = g^{n_1}\circ f^{n_1}\bigl(\psi(mP_0)\bigr) 
= g^{n_1}\circ f^{n_2}(mP_0) \in \Isog(X)\cdot P_0.
\]
Hence\footnote{We use that notation that for any subgroup~$Z\subset{X(\Kbar)}$, the division subgroup associated to~$Z$ is
$
Z^{\operatorname{div}} 
:= \bigl\{ P \in X(\Kbar): \text{$nP\in Z$ for some integer $n\ge1$} \bigr\}. 
$
For example, $\{0\}^{\operatorname{div}}=X(\Kbar)_\tors$.}
\[
\psi(mP_0) \in \Bigl( \Isog(X)\cdot P_0 \Bigr)^{\operatorname{div}}
\quad\text{for all $m\in\Mcal$.}
\]
To recapitulate, we have proven that
\begin{equation}
\label{eqn:YcapIsgoXP01}
\psi(\Mcal\cdot P_0)
\subseteq \Qcal \cap \Bigl( \Isog(X)\cdot P_0 \Bigr)^{\operatorname{div}}
\subseteq Y \cap \Bigl( \Isog(X)\cdot P_0 \Bigr)^{\operatorname{div}} .
\end{equation}
\par
The ring of isogenies~$\Isog(X)$ is a finitely generated~$\ZZ$-module. Indeed, it is a classical result that~$\Isog(X)$ is an order in a central simple algebra of rank at most~$2\dim(X)$; see for example~\cite[section~19]{MR2514037}. Hence the subgroup~$\Isog(X)\cdot{P}$ is a finitely generated subgroup of~$X$. We now apply~\cite{MR1323985}, which represents a culmination of fundamental work of Faltings, Vojta, Raynaud, Hindry, and others on the intersection of a subvariety of a semi-abelian variety with a subgroup of finite type. This result tells us that
\begin{equation}
\label{eqn:YcapIsgoXP02}
\overline{Y \cap \Bigl( \Isog(X)\cdot P_0 \Bigr)^{\operatorname{div}}}
= 
\left(\;\text{\parbox{.5\hsize}{the union of finitely many translates of abelian subvarieties of $X$}}\;\right).
\end{equation}
It follows from~\eqref{eqn:YcapIsgoXP01} that the set~\eqref{eqn:YcapIsgoXP02} contains~$\psi(\Mcal\cdot{P_0})$. We know that~$\psi$ is injective, and we proved earlier that~$\Mcal\cdot{P_0}$ is an infinite set. It follows that at least one of the abelian subvarieties appearing in the set~\eqref{eqn:YcapIsgoXP02} is positive dimensional. On the other hand, since~$Y$ is Zariski closed, it contains the set~\eqref{eqn:YcapIsgoXP02}, so we conclude that~$Y$ contains a translate of a positive dimensional abelian subvariety of~$X$. Our assumption that~$X$ is geometrically simple implies~$X$ has no non-trival positive dimensional abelian subvarieties. Hence~$Y=X$, which concludes our proof that~$f$ has propagation property~(C2$\forall$).
\par
We know that~$f$ has propagation property~(A) by assumption, we have proven that~$f$ has propagation property~(C2$\forall$), and we have proven~\eqref{eqn:Xfdenseisdenseabvar} that~$X_f^\dense(K)$ is Zariski dense in~$X$. It follows from Lemma~\ref{lemma:C2forallZDimpliesC3forall} that~$f$ has propagation property~(C3$\forall$).
\Case{2}{$\Ker(\f-1)=X$}
In this case~$\f=1$, so~$f$ is a pure translation
\[
f(x)=x+Q_0.
\]
Hence for all $n\in\NN$, we have~$f^n(x)=x+nQ_0$, which shows that
\begin{equation}
\label{eqn:OfxxplusNQ0}
\Orbit_f(x) = x + \NN\cdot Q_0.
\end{equation}
It follows that for any~$P\in{X}$ we have
\[
\Orbit_f(P) = P+\NN\cdot Q_0 = (P-P_0)+(P_0+\NN\cdot Q_0) = (P-P_0)+\Orbit_f(P_0),
\]
i.e., all~$f$-orbits are translates of one another. Hence the fact that~$\Orbit_f(P_0)$ is Zariski dense in~$X$ implies that every~$f$-orbit is Zariski dense in~$X$. This tells us that
\begin{equation}
\label{eqn:abvartransXfdenseX}    
\overline{X_f^\dense(K)} = \overline{X(K)} = X,
\end{equation}
where the second equality follows from the fact that~$\Orbit_f(P_0)$ is a Zariski dense subset of~$X(K)$.
\par
We next replace~$K$ with a finite extension, which by abuse of notation we again denote by~$K$, such that~$\rank{X(K)}\ge2$. This is always possible, since more generally it is known~\cite[Theorem~10.1]{MR337997} that~$\rank{X(\Kbar)}=\infty$. We let~$Q_1\in{X(K)}$ be a point such that~$Q_0$ and~$Q_1$ are~$\ZZ$-linearly independent. 
We claim that~$\NN\cdot{Q_1}$ has the following three properties:
\begin{parts}
\Part{(1)}
Every point in $\NN\cdot{Q_1}$ has a Zariski dense $f$-orbit.
\Part{(2)}
The points in $\NN\cdot{Q_1}$ generate distinct grand $f$-orbits.
\Part{(3)}
The set $\NN\cdot{Q_1}$ is Zariski dense in~$X$.
\end{parts}
The verification of~(1) is immediate, since as already noted, every point in~$X$ has Zariski dense $f$-orbit. For~(3), we note that since~$Q_1$ is a non-torsion point, any subvariety~$Y$ of~$X$ containing~$\NN\cdot{Q_1}$ would have infinite intersection with the finitely generated subgroup~$\ZZ\cdot{Q_1}$ of~$X$. Faltings's theorem~\cite{MR1109353} and the assumed simplicity of~$X$ then implies that~$Y=X$, which proves that~$\NN\cdot{Q_1}$ is Zariski dense.\footnote{We probably don't need the full strength of Faltings's theorem, since the assumption that~$Y$ contains all of~$\NN\cdot{Q_1}$ is much stronger than the assumption that $\#\bigl(Y\cap(\NN\cdot{Q_1})\bigr)=\infty$. It would suffice to show that the Zariski closure of $\NN\cdot{Q_1}$ is an algebraic subgroup of $X$.}
\par
It remains to prove~(2). For~$mQ_1,m'Q_1\in\NN\cdot{Q_1}$ we compute
\begin{align*}
\Orbit_f&(mQ_1)^\grand \cap \Orbit_f(m'Q_1)^\grand \ne \emptyset\\
&\quad\Longleftrightarrow\quad
f^n(mQ_1) = f^{n'}(m'Q_1) \quad\text{for some $n,n'\in\NN$,} \\
&\quad\Longleftrightarrow\quad
mQ_1+nQ_0 = m'Q_1+n'Q_0 \quad\text{for some $n,n'\in\NN$,} \\
&\quad\Longleftrightarrow\quad
(m-m')Q_1= (n'-n)Q_0 \quad\text{for some $n,n'\in\NN$,} \\
&\quad\Longrightarrow\quad
m=m'\quad\text{since $Q_0$ and $Q_1$ are $\ZZ$-linearly independent.}
\end{align*}
This completes the proof that distinct points in~$\NN\cdot{Q_1}$ have distinct grand~$f$-orbits.
\par
The properties~(1),~(2),~(3) that we proved for the set~$\NN\cdot{Q_1}$ are exactly that same as the properties~(1),~(2),~(3) that we proved for the set~$\Mcal\cdot{P_0}$ while proving Case~1. The remainder of the proof that~$f$ has propagation property~(C3$\forall$) in Case~2 is the same, \emph{mutatis mutandis}, as the proof in Case~1 with the set~$\Mcal\cdot{P_0}$ in Case~1 replaced with the set~$\NN\cdot{Q_1}$ for Case~2.
\end{proof}

\section{\texorpdfstring
{{\'E}tale Maps and $p$-adic Dimension}
{{\'E}tale Maps and p-adic Dimension}
} 
\label{section:etalepadic}
In this section we use $p$-adic arguments to prove orbit propagation results for {\'e}tale maps. 

\begin{definition}
We need to explain what we mean by the \emph{$p$-adic dimension} of an algebraic variety. For a smooth variety, we use the same analytic definition as for~$\RR$ and~$\CC$, namely we use local charts and locally convergent power series, and we declare that~$\AA^n(\QQ_p)$ has $p$-adic dimension~$n$. For a singular variety, we write the variety as a union of smooth varieties by inductively removing singular loci, and then the $p$-adic dimension of~$X(\QQ_p)$ is the maximum of the $p$-adic dimensions of the smooth varieties appearing in the decomposition. And since we will be using two notions of dimension in the section, we write~$\dimpadic$ for the $p$-adic dimension of a $p$-adic set, and we write~$\dimKrull$ for the usual Krull, or algebraic, dimension of an algebraic variety.
\end{definition}

\begin{lemma}
\label{lemma:padicdimen}
\begin{parts}
\Part{(a)}
Let~$X/\QQ_p$ be a smooth algebraic variety. Then
\[
X(\QQ_p)\ne\emptyset \quad\Longrightarrow\quad
\dimpadic X(\QQ_p) = \dimKrull X.
\]
\Part{(b)}
Let~$K$ be a number field, and let~$X/K$ be a smooth algebraic variety. Then for all but finitely many primes~$v$ of~$K$ satisfying~$K_v=\QQ_v$, we have
\[
\dimpadic X(K_v) = \dimKrull X.
\]
\Part{(c)}
Let~$X/\QQ_p$ be a smooth algebraic variety, and let~$Z\subsetneq{X}$ be a proper \textup(not necessarily smooth\textup) algebraic subset of~$Z$ defined over~$\QQ_p$. Then
\[
\dimpadic Z(\QQ_p) < \dimKrull X.
\]
\Part{(d)}
Let $K$ be a number field, let $X/K$ be a smooth projective variety with~$X(K)\ne\emptyset$, and suppose that there are infinitely many degree~$1$ primes~$v$ of~$K$ such that~$X(K)$ is $v$-adically dense in~$X(K_v)$. Then~$X(K)$ is Zariski dense in~$X$.
\end{parts}
\end{lemma}
\begin{proof}
(a)\enspace
Let~$n=\dimKrull X$ and~$Q\in{X}(\QQ_v)$. 
The smoothness of~$X$ and the $v$-adic implicit function theorem~\cite[page~73]{MR2179691} say that there is neighborhood of~$Q$ that is $v$-adic analytically
isomorphic to a neighborhood of the origin in~$\AA^n(\QQ_v)$. 
\par\noindent(b)\enspace
We choose a model~$\Xcal$ for~$X$ over the ring of integers of~$K$.
The smoothness of~$X$ tells us that the reduction of~$\Xcal$ modulo~$v$
is smooth for all but finitely many primes~$v$. If $\operatorname{char}(v)$ is sufficiently large, then the standard Lang--Weil estimate~\cite{MR65218} shows that~$X(\FF_v)\ne\emptyset$, and then smoothness and Hensel's lemma show that~$X(K_v)\ne\emptyset$. Restricting to primes with~$K_v=\QQ_v$, the desired result follows from~(a).
\par\noindent(c)\enspace  
We can check the assertion on Zariski open sets, and by Noetherian induction
on the dimension, we are reduced to proving that if~$V/\QQ_p$ is a smooth
(affine) variety, then
\[
\dimpadic V(\QQ_p) \le \dimKrull V.
\]
If~$V(\QQ_p)=\emptyset$, this is vacuously true. Otherwise we may apply~(a).
\par\noindent(d)\enspace  
We write~$\Closure_v$ for the $v$-adic closure of a set and~$\Closure_Z$ for the Zariski closure. So our assumption is that there are infinitely many degree~$1$ primes~$v$ with~$\Closure_v\bigl(X(K)\bigr)=X(K_v)$. We use this assumption and~(b) to deduce that there is a degree~$1$ prime~$v$, say of characteristic~$p\ge3$, such that
\begin{equation}
\label{eqn:dpCXKdpXKvdKX1}
\dimpadic \Closure_v\bigl(X(K)\bigr)
= \dimpadic X(K_v) = \dimKrull X.
\end{equation}
On the other hand, if $\Closure_Z\bigl(X(K)\bigr)$ is a proper subset of~$X$, then~(c) says that
\begin{equation}
\label{eqn:dpCXKdpXKvdKX2}
\dimpadic \bigl(\Closure_Z(X(K))\bigr)(K_v) < \dimKrull X
\quad\text{(strict inequality).}
\end{equation}
Since~\eqref{eqn:dpCXKdpXKvdKX1} and~\eqref{eqn:dpCXKdpXKvdKX2} are incompatible, we conclude that~$\Closure_Z\bigl(X(K)\bigr)$ is equal to~$X$.
\end{proof}

\begin{proposition} 
\label{proposition:schemeoverZp}
Let~$p\ge 3$ be prime, let~$\pi:\Xcal\to\Spec \ZZ_p$ be a smooth projective morphism with geometrically irreducible fibres, and let~$F:\Xcal\to \Xcal$ be a finite {\'e}tale~$\ZZ_p$-morphism. Let~$X=\Xcal\times_{\ZZ_p}\QQ_p$ be the generic fibre of~$\pi$, and let~$f:X\to X$ be the restriction of~$F$ to the generic fiber. Let~$x\in X(\QQ_p)$. Then the~$p$-adic closure of~$\Orbit_f(x)$ in~$X(\QQ_p)$ consists of finitely many~$p$-adic arcs, so in particular it has~$p$-adic dimension~$1$.
\end{proposition}
\begin{proof} 
We let~$X'=\Xcal\times_{\ZZ_p}\FF_p$ be the special fiber of~$\pi$, and 
we let~$n$ be the relative dimension $n$ of~$\pi$, so~$X$ and~$X'$ are smooth projective varieties of dimension~$n$. By properness,  the point $x\in{X(\QQ_p)}$ extends to a section $\s_x:\Spec\ZZ_p\to\Xcal$ of $\pi$, and the reduction of~$\s_x$ to the special fiber gives a point $x'\in X'(\F_p)$. Since $X(\F_p)$ is finite, there is no loss of generality if we assume that $f(x')=x'$.
\par
Expanding $f$ locally on the residue disk of $x'$ centered along $\sigma_x$, we see that $f$ is given by an $n$-tuple of convergent power series $g_1,...,g_n\in \ZZ_p[\![t_1,...,t_n]\!]$, where the $t_j$ are local parameters along $\sigma$. As in~\cite[Proposition~2.2]{MR2766180}, we may choose the $t_j$ so that the $g_j$ are congruent to linear polynomials modulo $p$. Since $f$ is \'etale, the Jacobian of the $g_j$ reduces modulo $p$ to a matrix in $\GL_n(\FF_p)$ that is necessarily of finite order. Replacing $f$ by a suitable iterate, we get that 
\[
(g_1,...,g_n) \equiv (t_1,...,t_n) \pmod{p}.
\]
It now follows from~\cite{MR3210707} that the map
\[
\NN\longrightarrow X(\QQ_p),\quad
n\longmapsto f^n(x),
\]
extends to a locally $p$-adic analytic map
\begin{equation}
\label{eqn:ZptoXQp}
\ZZ_p\longrightarrow X(\QQ_p). 
\end{equation}
Hence the $p$-adic closure of~$\Orbit_f(x)$ is contained in the image of~\eqref{eqn:ZptoXQp}.
\end{proof}

\begin{theorem} 
\label{theorem:XQpdim2}
Let $K$ be a number field, let $X/K$ be a smooth projective variety of dimension at least~$2$, and let $f:X\to X$ be a finite {\'e}tale morphism. Suppose that there are infinitely many degree~$1$ primes $v$ of $K$, i.e., primes satisfying $K_v=\QQ_v$, such that~$X(K)$ is $v$-adically dense in~$X(\QQ_v)$. Then~\textup{(B1)} in Table~$\ref{table:propogationstatements}$ is true for~$(X,f)$, i.e., if~$\G_1,\ldots,\G_r\subset{X(K)}$ are~$f$-orbits, then
\[
\text{$X(K)\setminus(\G_1\cup\cdots\cup\G_r)$
is Zariski dense in $X$.}
\]
\end{theorem}

Before proving Theorem~\ref{theorem:XQpdim2}, we state and prove two corollaries, one of which will require a preliminary lemma.

\begin{corollary}
\label{corollary:CoroRat}
Let~$X/K$ be a smooth rational variety defined over a number field, and let~$f:X\to{X}$ be a finite {\'e}tale morphism defined over~$K$. Then~\textup{(B1)} in Table~$\ref{table:propogationstatements}$ is true for~$(X,f)$. In particular,~\textup{(B1)} is true for automorphisms of rational varieties.
\end{corollary}

\begin{corollary}
\label{corollary:B1foretaleabelianquotients}
Let~$X$ be an {\'e}tale quotient of an abelian variety,\footnote{In other words, there is an abelian variety~$A$ and a finite group of automorphisms~$G\subset\Aut(A)$ (not necessarily isogenies) so that $X\cong{A}/G$ and $A\to{X}$ is {\'e}tale.} and let $f:X\to{X}$ a morphism. Then~\textup{(B1)} in Table~$\ref{table:propogationstatements}$ is true for~$(X,f)$. In particular,~\textup{(B1)} is true for abelian varieties and for bielliptic surfaces.
\end{corollary}

\begin{proof}[Proof of Corollary~$\ref{corollary:CoroRat}$]
Possibly after replacing~$K$ with a finite extension, we can find a birational map $\f:\PP^N\dashrightarrow{X}$ defined over~$K$, and we let~$U\subseteq\PP^N$ be a non-empty Zariski open subset defined over~$K$ such that~$\f|_U$ is an isomorphism onto its image. The fact that~$U$ is an open subset of~$\PP^N$ implies that~$U(K)$ is $v$-adically dense in~$U(K_v)$ for all places of~$K$. (This follows simply from the fact that~$K$ is $v$-adically dense in~$K_v$.) It follows that~$\f\bigl(U(K)\bigr)=\f(U)(K)$ is~$v$-adically dense in~$\f(U)(K_v)$,
and then since~$\f(U)$ is a Zariski open subset of the irreducible variety~$X$, we see that~$\f(U)(K)$ is $v$-adically dense in~$X(K_v)$. Hence the large set~$X(K)$ is $v$-adically dense in~$X(K_v)$, and the desired conclusion follows from Theorem~\ref{theorem:XQpdim2}.
\end{proof}

\begin{lemma}
\label{lemma:etalequotabelianvariety}
Let~$A$ be abelian variety, let $\pi:A\to B$ be an {\'e}tale map, and let $f:B\to B$ be a surjective morphism. Then $f$ is {\'e}tale.  
\end{lemma}
\begin{proof}
We let~$X$ be the fiber product of~$A$ and~$B$ relative to the map~$\pi$ and~$f$, so we have a commutative diagram
\[
\begin{CD}
  X @>h>> A \\
  @V p VV @V \pi VV \\
  B @>f>> B \\
\end{CD}
\]
The assumption that~$\pi$ is {\'e}tale implies that~$X$ is smooth and~$p$ is {\'e}tale.
\par
The fact that~$p$ is {\'e}tale tells us that~$\kappa(X)=0$ and that~$\Omega^1_X$ is semiample, since~$B$ has these properties, so~\cite{MR1183561} allows us to conclude that~$X$ is the quotient of an abelian variety~$C$. Considering the composite ma $C\to{X}\to{A}$, we conclude that~$X\to{A}$ is {\'e}tale, and thus that $f:B\to{B}$ is also {\'e}tale.
\end{proof}

\begin{proof}[Proof of Corollary $\ref{corollary:B1foretaleabelianquotients}$]
Lemma~\ref{lemma:etalequotabelianvariety} says that it is enough to check the case that~$X$ is an abelian variety. Since every map between dominant map~$X\to{X}$ of an abelian variety is {\'e}tale, Theorem~\ref{theorem:XQpdim2} says that it suffices to check that for all but finitely many places~$v$ of~$K$, the set of~$K$-rational points~$X(K)$ is $v$-adically dense in~$X(K_v)$. This will be true provided that the rank of the Mordell-Weil group of each simple isogeny factor of~$X(K)$ is sufficiently large with respect to the dimension of~$X$; see~\cite{MR2793038}. And we can make these ranks sufficiently large by taking a finite extension of~$K$.
\end{proof}

\begin{proof}[Proof of Theorem~$\ref{theorem:XQpdim2}$]
We choose a model~$\Xcal\to\Spec(\Ocal_K)$ of~$X$ over the ring of integers~$\Ocal_K$ of~$K$, and then~$f$ extends to a rational map~$F:\Xcal\dashrightarrow\Xcal$. For any finite set of places~$S$ of~$K$ that includes the archimedean places, we denote the ring of~$S$-integers by~$\Ocal_{K,S}$. We choose such a set~$S$ having the property that the scheme~$\Xcal_S=\Xcal\times_{\Ocal_K}\Ocal_{K,S}$ is smooth and proper over~$\Ocal_{K,S}$ with geometrically irreducible fibers and such that the map~$F_S:\Xcal_S\to\Xcal_S$ is a dominant {\'e}tale morphism. Since the degree~$1$ primes have positive density, the assumption of Theorem~\ref{theorem:XQpdim2} and the finiteness of~$S$ imply that we can find a degree~$1$ prime~$v\notin{S}$ of characteristic $p\ge3$ such that $X(K)$ is $p$-adically dense in $X(\QQ_p)$. Using the smoothness of~$X$, the Lang--Weil estimate~\cite{MR65218}, and Hensel's lemma, we may also assume that~$X(\QQ_p)\ne\emptyset$. 
\par
Using the closure notation from the proof of Lemma~\ref{lemma:padicdimen}(d), we have
\begin{align}
\label{eqn:2ledimKrullX1}
2 &\le \dimKrull X
&&\text{by assumption,} \notag\\
&= \dimpadic X(\QQ_p)
&&\text{from Lemma~\ref{lemma:padicdimen}(a),} \\
\label{eqn:2ledimKrullX2}
&= \dimpadic \Closure_p\bigl(X(K)\bigr)
&&\begin{tabular}[t]{@{}l}
since $\Closure_p\bigl(X(K)\bigr)=X(\QQ_p)$\\ by assumption.\\
\end{tabular}
\end{align}
\par
On the other hand, for~$x\in{X(K)}$, it follows from Proposition~\ref{proposition:schemeoverZp} that 
\[
\dimpadic \Closure_p\bigl( \Orbit_f(x) \bigr) \le 1.
\]
To ease notation, we let
\[
\G := \G_1\cup\cdots\cup\G_r
\]
be the finite union of orbits in the statement of the Corollary. Then, since the dimension of a finite union of sets is the maximum of the dimensions of the components,  we have
\begin{equation}
\label{eqn:QpdimGammai}
\dimpadic \Closure_p(\G) \le 1.
\end{equation}
Thus~\eqref{eqn:2ledimKrullX1} says that $\Closure_p\bigl(X(K)\bigr)$ has $p$-adic dimension at least~$2$, while~\eqref{eqn:QpdimGammai} says that~$\Closure_p(\G))$ has $p$-adic dimension at most~$1$, so the complement satisfies
\[
\dimpadic\Bigl(\Closure_p\bigl(X(K)\bigr) \setminus \Closure_p(\G)) \Bigr)
= \dimpadic \Closure_p\bigl(X(K)\bigr).
\]
Applying~\eqref{eqn:2ledimKrullX2} gives
\begin{equation}
\label{eqn:dimpCpXKnegCpGdimK}
\dimpadic\Bigl(\Closure_p\bigl(X(K)\bigr) 
\setminus \Closure_p(\G) \Bigr)
= \dimKrull X.
\end{equation}
Elementary topology\footnote{It is an exercise to check that if~$U$ and~$V$ are subsets of a topological space, then their closures satisfy
$\Closure(U\setminus V)\supseteq\Closure(U)\setminus\Closure(V)$. See for example~\cite{MSE493989}.} tells us that
\[
\Closure_p\Bigl(X(K) \setminus \G\Bigr)
\supseteq
\Closure_p\bigl(X(K)\bigr) 
\setminus \Closure_p\bigl(\G\bigr),
\]
so taking $p$-adic dimension and using~\eqref{eqn:dimpCpXKnegCpGdimK} yields
\begin{equation}
\label{eqn:dimCpXKnegGgedimKX}
\dimpadic\Bigl( \Closure_p\bigl(X(K) \setminus \G\bigr) \Bigr) \ge \dimKrull X.
\end{equation}
\par
Let
\[
Z := \Closure_Z\bigl( X(K)\setminus\G \bigr). 
\]
We note that
\begin{align*}
\dimpadic & Z(\QQ_p)  \\
&\ge \dimpadic \Closure_p\bigl( X(K)\setminus\G \bigr)
\quad\begin{tabular}[t]{@{}l}
since $Z$ is a Zariski closed\\ set that contains $X(K)\setminus\G$,\\
\end{tabular}
\\
&\ge \dimKrull X
\quad\text{from \eqref{eqn:dimCpXKnegGgedimKX}.}
\end{align*}
On the other hand, if~$Z\subsetneq{X}$ is a proper subset of~$X$, then Lemma~\ref{lemma:padicdimen}(c) says that
\[
\dimpadic Z(\QQ_p) < \dimKrull X.
\]
Hence~$Z=X$, which is the desired conclusion.
\end{proof}

\section{Surfaces}
\label{section:surfaces}

Our goal in this section is to prove an orbit propagation result for all smooth projective surfaces. We note that we have already proven a number of cases that apply in particular to surfaces. These results are summarized in Table~\ref{table:surfaceorbitpropagation}.

\begin{table}[ht]
\begin{center}
\begin{tabular}{|lc@{\quad}l|} \hline
$\PP^2$, $\deg(f)\ge2$ 
& $\text{(A)}\Longrightarrow\text{(C1) and (B1$\exists$)}$
& Theorem~\ref{theorem:projectivespacedeg2}(a,b) \\ \hline
\begin{tabular}{@{}l} $\PP^2$, $\deg(f)=1$\\  \end{tabular}
& $\text{(A)}\Longrightarrow\text{(C3$\forall$)}$ 
& Theorem~\ref{theorem:projectivespacedeg2}(c)\\ \hline
\begin{tabular}{@{}l} Geometrically simple\\ abelian surfaces\\ \end{tabular}
& $\text{(A)}\Longrightarrow\text{(C3$\forall$)}$
& Theorem~\ref{theorem:simpleabelianvariety} \\ \hline
K3 surfaces 
& $\text{(A)}\Longrightarrow\text{(C1)}$ 
& Theorem~\ref{theorem:K3} \\ \hline
\begin{tabular}{@{}l} Rational surface\\ $f$ \'etale\\ \end{tabular}
& $\text{(A)}\Longrightarrow\text{(B1)}$ 
& Corollary~\ref{corollary:CoroRat} \\ \hline
\begin{tabular}{@{}l} \'etale quotient of\\ an abelian surface\\ \end{tabular}
& $\text{(A)}\Longrightarrow\text{(B1)}$ 
& Corollary~\ref{corollary:B1foretaleabelianquotients} \\ \hline
\end{tabular}
\end{center}
\caption{Orbit propagation results that apply to surfaces} 
\label{table:surfaceorbitpropagation}
\end{table}

\begin{theorem}
\label{theorem:surfaces}
Let~$K$ be a number field and let~$X/K$ be a smooth projective surface. Then possibly after replacing~$K$ by a finite extension, for every finite collection of $f$-orbits $\G_1,\ldots,\G_r\subset{X}(K)$, we have
\[
\overline{X(K) \setminus (\G_1\cup\cdots\cup\G_r)} = X.
\]
In the terminology of Table~\textup{\ref{table:propogationstatements}}, smooth projective surfaces satisfy the orbit propagation statement
\[
\textup{(A)} \quad\Longrightarrow\quad \textup{(B1)}.
\]
\end{theorem}

Before starting the proof of Theorem~\ref{theorem:surfaces}, we
give some lemmas that will be used in the proof.

\begin{lemma} 
\label{lemma:Lemmakmi} 
Let $Y$ be a relatively minimal surface with $\kappa(Y)=-\infty$. Let $X$ be obtained from $Y$ by a sequence of blow-ups. If $f:X\to X$ is a non-invertible surjective morphism, then some iterate of $f$ is induced from an endomorphism of $Y$. In particular, if $\kappa(X)=-\infty$, then Theorem~$\ref{theorem:surfaces}$ for all relatively minimal surfaces $Y$ and and all non-invertible maps implies Theorem~$\ref{theorem:surfaces}$ for $X$ with a non-invertible map $f$.
\end{lemma}
\begin{proof}
See \cite[Lemma~4.1]{MR3871505}.
\end{proof}

\begin{lemma} 
\label{lemma:negonecurves}
Let $X$ be a smooth projective complex surface that is not birational to $\PP^2$. Then $X$ contains at most finitely many $(-1)$-curves.
\end{lemma}
\begin{proof}
This is well-known. See for example~\cite{MathOverflow267708}.
\end{proof}

We next consider a version of Lemma~\ref{lemma:Lemmakmi} for non-rational surfaces.

\begin{lemma}
\label{lemma:LemmaNonMin} 
Let $Y$ be a minimal surface that is not birational to $\PP^2$. Let $X$ be obtained from $Y$ by a sequence of blow-ups. If $f:X\to X$ is a  surjective morphism, then some iterate of $f$ is induced from an endomorphism of $Y$. 
In particular, Theorem~$\ref{theorem:surfaces}$ for $Y$ implies Theorem~$\ref{theorem:surfaces}$ for $X$.
\end{lemma}
\begin{proof} 
If~$f$ is an automorphism, then it maps exceptional curves into exceptional curves and we are done. So we assume that~$f$ is non-invertible and consider two cases. First, if $\kappa(X)\ge 0$, then~\cite[Lemma~4.2]{MR3871505} tells us that $X=Y$. Second, if $\kappa(X)=-\infty$, then Lemma~\ref{lemma:Lemmakmi} gives the desired result.
\end{proof}

\begin{proof}[Proof of Theorem~$\ref{theorem:surfaces}$]
We note from Lemma~\ref{lemma:LemmaNonMin} that Theorem~\ref{theorem:surfaces} for minimal surfaces implies Theorem~\ref{theorem:surfaces} for all surfaces, except possibly in the case that~$X$ is a rational surface and $f:X\to X$ is an automorphism. But that case is covered by Corollary~\ref{corollary:CoroRat}. So we have reduced the proof of Theorem~\ref{theorem:surfaces} to the case of minimal surfaces. The proof of Theorem~\ref{theorem:surfaces} is now a case-by-case analysis via the classical classification of surfaces~\cite[Section~V.6]{MR0463157}.
\par
We remark that all of the implications in Table~\ref{table:surfaceorbitpropagation} are at least as strong as the assertion of Theorem~\ref{theorem:surfaces}, since each of the properties~(C1),~(C3$\exists$), and~(C3$\forall$) imply~(B1); cf.\ Table~\ref{table:propagationimplications}. Hence Theorem~\ref{theorem:surfaces} is true for the surfaces listed in Table~\ref{table:surfaceorbitpropagation}.
\par
\Case{1}{$\kappa(X)=-\infty$, $X$ is rational}
Since we are assuming that~$X$ is minimal, this implies that~$X\cong\PP^2$, so Theorem~\ref{theorem:projectivespacedeg2}(a,b,c) give something stronger than the desired result.
\Case{2}{$\kappa(X)=-\infty$, $X$ is ruled}
Let~$\pi:X\to{B}$ be a surjective map from~$X$ to a curve~$B$ whose fibers are~$\Kbar$-isomorphic to~$\PP^1$. If~$g(B)\ge2$, then~$B(K)$ is finite, so there are no dense orbits. We may thus assume that~$g(B)$ is~$0$ or~$1$.
\par
Possibly after replacing $f$ by $f^2$ (cf.\ \cite[Lemma~5.4]{MR3871505}), we may assume that $f$ is semi-conjugate to $\pi$, i.e, there is a map $h:B\to B$ and a commutative diagram
\begin{equation}
\label{eqn:ruledsurface}
\begin{CD}
X @>f>> X \\
@V \pi VV @V \pi VV \\
B @>h>> B \\
\end{CD}
\end{equation}
Furthermore, $\pi$ has a section $\sigma:B\to X$; see \cite[Lemma~5.1]{MR3871505}. Replacing $K$ by a finite extension,  we may assume that everything is defined over $K$. Then for every finite extension~$L/K$, we have
\[
\pi\bigl(X(L)\bigr) \supseteq \pi\bigl(\s\bigl(B(L)\bigr)\bigr) = B(L).
\]
Thus$\pi\bigl(X(L)\bigr)=B(L)$, where~$B$ is either~$\PP^1$ or an elliptic curve. The implication (A)$\Rightarrow$(B1) is true for the curve~$B$ by the dimension~$1$ cases of Theorem~\ref{theorem:projectivespacedeg2} and~\ref{theorem:simpleabelianvariety}, where for the latter we note that an elliptic curve is always geometrically simple.
\par
We want to prove that~(A) implies~(B1) for~$X$, so we assume that there is a point~$x_0\in{X(K)}$ whose orbit~$\Orbit_f(x_0)$ is Zariski dense in~$X$. It follows from the commutative diagram~\eqref{eqn:ruledsurface} that~$\Orbit_h\bigl(\pi(x_0)\bigr)$ is Zariski dense in~$B$, so as noted in the previous paragraph, we know that~(B1) is true for~$B$.
\par
Let~$\G_1,\ldots,\G_r$ be $f$-orbits in~$X(K)$. Then~$\pi(\G_1),\ldots,\pi(\G_r)$ are~$h$ orbits in~$B$, so the validity of~(B1) for~$B$ tells us that
\[
B^\circ(K) := B(K) \setminus \bigl( \pi(\G_1),\ldots,\pi(\G_r) \bigr)
\quad\text{is Zariski dense in $B$.}
\]
The fibers of~$\pi$ over the points in~$B^\circ(K)$ are $K$-isomorphism to~$\PP^1$, since the fiber over~$b\in{B^\circ(K)}$ is a rational curve $X_b:=\pi^{-1}(b)$ that contains the $K$-rational point~$\s(b)$. Hence~$X_b(K)$ is Zariski dense in~$X$, so the density of~$B^\circ(K)$ in~$B$ implies that
\[
\bigcup_{b\in B^\circ(K)} X_b(K)\quad\text{is Zariski dense in $X$.}
\]
This set is disjoint from~$\G_1\cup\cdots\cup\G_r$ be construction, which completes the proof that~$X(K)\setminus(\G_1\cup\cdots\cup\G_r)$ is Zariski dense in~$X$, and thus the proof that~(B1) holds for~$X$.
\Case{3}{$\kappa(X)=0$, $X$ is a K3 surface}
Theorem~\ref{theorem:K3} proves  $\text{(A)}\Longrightarrow\text{(C1)}$ for~K3 surfaces, which is stronger than the desired result.
\Case{4}{$\kappa(X)=0$, $X$ is an Enriques surface}
There is an \'etale quotient map $\pi:X'\to X$ with $X'$ a $K3$ surface, and the endomorphism $f:X\to X$ lifts to $X'$ because~K3 surfaces are simply connected. 
Hence the desired result for~$X$ follows from the fact that it is true the~K3 surface~$X'$.
\Case{5}{$\kappa(X)=0$, $X$ is an abelian surface}
This is covered by Corollary~\ref{corollary:B1foretaleabelianquotients}.
\Case{6}{$\kappa(X)=0$, $X$ is a bielliptic surface}
This is covered by Corollary~\ref{corollary:B1foretaleabelianquotients}.
\Case{7}{\begin{tabular}[t]{@{}l@{}}
$\kappa(X)=1$, $X$ is an elliptic surface\\
$\kappa(X)=2$, $X$ is of general type.\\
\end{tabular}}
In general, if~$\kappa(X)>0$, then we can use Nakayama and Zhang~\cite[Theorem~A]{MR2551469}. They prove that any dominant rational self-map of a smooth projective variety of positive Kodaira dimension factors over a positive dimensional base, and hence does not have any Zariski dense orbits. See also~\cite[Corollary~8.2]{MR3871505} for a more elementary proof of the case that we need. Thus all versions of orbit propagation are vacuously true, since there are no Zariski dense orbits. (We mention that varieties of general type cannot even have endomorphisms of infinite order, so orbits are necessarily finite in that case, making orbit propagation even more vacuous.)
\end{proof}


\appendix

\section{The Zariski Dense Orbit Conjecture}
\label{section:ZDconjecture}

In this section, for the convenience of the reader, we discuss various versions of the Zariski Density Conjecture. As noted in the introduction and in Figure~\ref{figure:ZDCimpliesOPC}, these conjectures are complementary to our Orbit Propagation Conjectures, since the former says that a geometric property implies the existence of Zariski dense orbits over an algebraically closed field, while the latter says that the existence of one such orbit implies that there are many Zariski dense orbits defined over a finitely generated field. To do full justice to the Zariski Density Conjecture, we work in somewhat greater generality than in the other parts of this paper, using the notation described in Table~\ref{table:ZDnotation}.

\begin{table}[ht]
\begin{center}
\framebox{\parbox{.95\hsize}{
\begin{itemize}
\setlength{\itemsep}{5pt}
\item[$\bfk$\enspace]
an algebraically closed field of characteristic~$0$ with\\
\text{$\operatorname{trdeg_\QQ(\bfk)<\infty}$} (for example $\bfk=\Qbar$)
\item[$X/\bfk$\enspace]
an irreducible (quasi-projective) variety defined over $\bfk$
\item[$f$\enspace]
a dominant rational map $f:X\dashrightarrow{X}$ defined over~$\bfk$
\item[$\bfk(X)$\enspace]
the function field of $X$
\item[$\bfk(X)^f$]
the subfield of $f$-invariant functions, i.e., the subfield\\
\hspace*{4em}$\bfk(X)^f=\{\f\in\bfk(X):\f\circ f=\f\}$
\end{itemize}
}}
\end{center}
\caption{Notation for the Zariski density conjecture}
\label{table:ZDnotation}
\end{table}

\begin{definition}
\label{definition:fiberedmap}
The map~$f:X\dashrightarrow{X}$ is \emph{fibered} if there is a dominant rational map~$\f:X\to\PP^1$ satisfying~\text{$\f\circ{f}=\f$}; or equivalently, if~\text{$\bfk(X)^f\ne\bfk$}.
\end{definition}

Clearly the existence of a point~$x\in{X(\bfk)}$ whose $f$-orbit is well-defined and Zariski dense implies that~$f$ cannot be fibered. Building on a conjecture of Zhang~\cite{MR2408228},
Amerik--Bogomolov--Rovinsky and Medvedev--Scanlon independently suggested that the converse should  hold, and Xie formulated a stronger version.

\begin{conjecture}
\label{conjecture:ZDconjecture}
With notation as in Table~\ref{table:ZDnotation}\textup:
\begin{parts}
\Part{(a)}
\textup{(The Zariski-Density (ZD) Orbit Conjecture \cite[Conjecture~1.2]{MR2862064} and \cite[Conjecture 5.10]{MR3126567})}  
\[
X_f^\dense(\bfk)=\emptyset\quad\Longleftrightarrow\quad \text{$f$ is fibered}.
\]
\Part{(b)}
\textup{(The Strong Zariski-Density (SZD) Orbit Conjecture \cite[Conjecture~1.4]{Xie2022AdelicZariskiDensity})}  
Assume that the dominant rational map~$f$ is not fibered. Then for every non-empty Zariski open subset~${U}\subset{X}$, the set\footnote{Xie further conjectures that there is a point  \text{$x\in X_f^\dense(\bfk)\cap{U}$} satisfying~\text{$\Orbit_f(x)\subset{U}$.}}
\[
\text{$X_f^\dense(\bfk)\cap U$ is Zariski dense in $U$.}
\]
\end{parts}
\end{conjecture}

\begin{remark}
Xie~\cite[Section~1.3]{Xie2022AdelicZariskiDensity} notes that ``if we have one Zariski dense orbit, we expect many such orbits. However, the Zariski topology is too weak to describe such phenomena; in particular, one cannot expect to have a nonempty Zariski open set of points having Zariski dense orbits.''  He goes on to develop a very interesting \emph{adelic topology} on~$X(\bfk)$.\footnote{We note that Xie's adelic topology is not the restriction to~$X(\bfk)$ of the natural topology on the adelic points of~$X$. In particular, Xie's adelic topology is~$T_1$ (Fr{\'e}chet), but not~$T_2$ (Hausdorff), while the restriction topology is~$T_2$.} Xie's adelic topology has many agreeable properties~\cite[Section~1.2]{Xie2022AdelicZariskiDensity}. In particular, it is stronger than the Zariski topology, and thus has more
open sets, which allows him to make the following conjecture.
\end{remark}

\begin{conjecture}
\label{conjecture:SAZD}
\textup{(The Adelic Zariski-Density (AZD) Orbit Conjecture \cite[Conjecture~1.10]{MR2862064})}
With notation as in Table~\ref{table:ZDnotation}, if the dominant rational map~$f$ is not fibered, then there exists a non-empty adelic open subset~${U}\subseteq{X}(\bfk)$ such that
\[
U \subseteq X_f^\dense(\bfk).
\]
In other words, there is an adelic open set~$U$  such that every point in~$U$ has a well-defined Zariski dense forward orbit.
\end{conjecture}

\begin{remark}
\label{remark:adelicopen}
We thank Junyi Xie for the following remarks (private communication):
Let~$U\subset{X(\bfk)}$ be a non-empty adelic open set as in Conjecture~\ref{conjecture:SAZD}. Then~$U$ is automatically Zariski dense in~$X$, but if~$C$ is defined over a finitely generated subfield of~$\bfk$, it need not be true that there is a finitely generated field of definition~$K\subset\bfk$ of~$X$ such that~$U\cap{X(K)}$ is Zariski dense in~$X$. For example, if~$X/\QQ$ is a field of genus at least~$2$, any non-empty adelic set~$U\subset{X(\Qbar)}$ is Zariski dense in~$X$, but~$X(K)$ is finite for every number field~$K/\QQ$.
\par
Xie notes that it might be possible to modify the definition of the adelic topology so that Conjecture~\ref{conjecture:SAZD} for~$\PP^N$ is reasonable with the added condition that $U\cap\PP^N(L)$ is Zariski dense for some finitely generated subfield~$L$ of~$\bfk$. However, this alternative topology would probably not be open under flat morphisms.
\end{remark}

We describe some of the cases where the~ZD conjecture and its various generalizations have been proven. 

\begin{theorem}
  The ZD conjecture is true for the following classes of varieties and maps, and in many cases, the SZD conjecture and/or the AZD conjecture are true\textup{:}
\begin{parts}
\Part{(a)}
\cite[Corollary~9]{MR2784670}
The field~$\bfk$ is uncountable, e.g., $\bfk=\CC$ or $\CC_p$. 
\Part{(b)}
\cite[Theorem~1.11]{MR4530188},
\cite[Theorem 1.15]{Xie2022AdelicZariskiDensity}
$X/\bfk$ is an irreducible (non-singular) projective surface, and $f:X\to{X}$ is
a dominant endomorphism.
\Part{(c)}
\cite{MR4396628,MR3937327,MR3557780,Xie2022AdelicZariskiDensity}
$X/\bfk$ is a (semi)-abelian variety, and $f:X\to{X}$ is a dominant endomorphism.
\Part{(d)}
\cite[Theorem 1.13]{Xie2022AdelicZariskiDensity}
$X=\AA^N$ and $f:\AA^N\to\AA^N$ is a dominant endomorphism
of the form~$f(x_1,\ldots,x_N)=\bigl(f_1(x_1),\ldots,f_N(x_N)\bigr)$.
\Part{(e)}
\cite[Theorem 1.16]{Xie2022AdelicZariskiDensity} \textup(joint with T.\ Tucker\textup)
$X=(\PP^1)^N$ and $f:X\to{X}$ is a dominant endomorphism.
\Part{(f)}
\cite{MR3705271}
$X=\AA^2$ and $f:\AA^2\to\AA^2$ is a dominant (regular) endomorphism.
\Part{(g)}
\cite{MR3829742}
$X$ is a  connected commutative linear algebraic groups and $f:X\to{X}$ is a group endomorphism.
\end{parts}
See~\cite[Remark~1.13]{MR4530188} and~\cite[Section~11]{mengzhang2023} for further details.
\end{theorem}

\section{Additional Results}
\label{appendix:elementaryresults}

This section contains additional material that  illuminates or is used in the text.

\subsection{Elementary Properties of Orbits}
\label{appendix:elempropsorbits}

The following lemma describes some elementary properties of orbits.

\begin{lemma}
\label{lemma:orbitelemproperties}
\begin{parts}
\Part{(a)}
Grand $f$-orbit equivalence is an equivalence relation.
\Part{(b)}
Let~$\G_1$ and~$\G_2$ be~$f$-grand orbits. Then either
\[
\G_1=\G_2 \quad\text{or}\quad \G_1\cap\G_2=\emptyset.
\]
Hence for $P,Q\in X$, we have
\[
P\fgeq{Q} \quad\Longleftrightarrow\quad \Orbit_f^\grand(P)=\Orbit_f^\grand(Q).
\]
\Part{(c)}
If~$\G$ is a grand orbit, then either all of its points are preperiodic, or all of its points are wandering. 
\Part{(d)}
Let~$\G$ be a grand orbit. If~$\G$ contains a point with Zariski dense forward orbit, then every point in~$\G$ has Zariski dense forward orbit. 
\Part{(e)}
If~$f$ is an automorphism, then the $f$-grand orbit of~$Q$ is the union of orbits of~$f$ and it's inverse,
\[
\Orbit_f(Q)^\grand = \Orbit_f(Q) \cup \Orbit_{f^{-1}}(Q).
\]
\Part{(f)}
If~$f$ is an automorphism, then~$P$ is~$f$-wandering if and only if~$P$ is~$f^{-1}$-wandering.
\Part{(g)}
Let $P,Q\in{X}$. Then
\[
P\fgeq Q \quad\Longrightarrow\quad 
\Bigl( P\in X^\dense \;\Longleftrightarrow\; Q\in X^\dense \Bigr).
\]
\end{parts}
\end{lemma}

\begin{proof}[Proof of Lemma \textup{\ref{lemma:orbitelemproperties}}]
(a)\enspace
We first note that
\[
\Orbit_f(P)\cap\Orbit_f(P)\supseteq\{P\}\ne\emptyset
\quad\Longrightarrow\quad P\fgeq P.
\]
Next, since the definitions of grand $f$-orbit equivalence is symmetric in~$P$ and~$Q$, we have
\[
P\fgeq Q \quad \Longrightarrow\quad Q\fgeq P.
\]
It remains to check transitivity. We compute
\begin{align*}
P & \fgeq Q \quad\text{and}\quad Q\fgeq R \\
&\quad\Longrightarrow\quad
\Orbit_f^\grand(P)=\Orbit_f^\grand(Q) 
\quad\text{and}\quad
\Orbit_f^\grand(Q)=\Orbit_f^\grand(R) \\
&\quad\Longrightarrow\quad
\Orbit_f^\grand(P)=\Orbit_f^\grand(R) \\
&\quad\Longrightarrow\quad
P\fgeq R.
\end{align*}
\par\noindent(b)\enspace
Let $\G_1=\Orbit_f^\grand(P_1)$ and  $\G_2=\Orbit_f^\grand(P_2)$, and assume that
$\G_1\cap\G_2\ne\emptyset$. Let $Q\in\G_1\cap\G_2$. This implies that there are iterates of~$f$ satisfying
\[
f^{n_1}(P_1) = f^{m_1}(Q)
\quad\text{and}\quad
f^{n_2}(P_2) = f^{m_2}(Q).
\]
Now let~$R\in\G_1$ be an arbitrary point in~$\G_1$. Our goal is to show that~$R\in\G_2$. The assumption that~$R\in\G_1\Orbit_f^\grand(P_1)$ tells us that there are iterates of~$f$ satisfying
\[
f^p(P_1) = f^q(R).
\]
This allows us to compute
\[
f^{m_2+n_1+q}(R) = f^{m_2+n_1+p}(P_1) = f^{m_2+m_1+p}(Q) = f^{n_2+m_1+p}(P_2.)
\]
This shows that~$\Orbit_f(R)\cap\Orbit_f(P_2)\ne\emptyset$, and hence by definition
$R\in\Orbit_f^\grand(P_2)=\G_2$. This completes the proof that~$\G_1\subseteq\G_2$, and the proof of the opposite inclusion follows by symmetry. Hence~$\G_1=\G_2$.
\par\noindent(c,d)\enspace
Let~$P,Q\in\G$. By definition of grand orbit, the forward orbits of~$P$ and~$Q$ eventually merge, so their set difference\footnote{We use the standard set theory notation that the set difference of sets~$A$ and~$B$ is $A\operatorname{\scriptstyle\triangle}{B}:=(A\cup{B})\setminus(A\cap{B})$, i.e., it is the set of points in exactly one of~$A$ and~$B$.}
\[
\Orbit_f(P) \operatorname{\scriptstyle\triangle} \Orbit_f(Q)
\quad\text{is a finite set.}
\]
This implies first that
\[
\#\Orbit_f(P) - \#\Orbit_f(Q)\quad\text{is finite,}
\]
which shows that~$P$ and~$Q$ are either both preperiodic or both wandering, which completes the proof of~(c). Second, it implies that
\[
\overline{\Orbit_f(P)} \setminus \overline{\Orbit_f(Q)} \quad\text{is finite,}
\]
and similarly with~$P$ and~$Q$ reversed. Hence~$\Orbit_f(P)$ is Zariski dense if and only if~$\Orbit_f(Q)$ is Zariski dense, which completes the proof of~(d).
\par\noindent(e)\enspace
We compute
\begin{align*}
P \in \Orbit_f(Q)^\grand
&\quad\Longleftrightarrow\quad
f^n(P) = f^m(Q) \quad\text{for some $n,m\in\NN$,} \\
&\quad\Longrightarrow\quad
\begin{cases}
P = f^{m-n}(Q) \in \Orbit_f(Q) &\text{if $m\ge n$,} \\
P = (f^{-1})^{n-m}(Q) \in \Orbit_{f^{-1}}(Q) &\text{if $n\ge m$.} \\
\end{cases} \\
&\quad\Longrightarrow\quad
P \in \Orbit_f(Q) \cup \Orbit_{f^{-1}}(Q).
\end{align*}
This gives the inclusion
\[
\Orbit_f(Q)^\grand \subseteq \Orbit_f(Q) \cup \Orbit_{f^{-1}}(Q).
\]
For the opposite inclusion, we note that the forward orbit of~$Q$ is clearly contained in the grand orbit of~$Q$. Finally, if~$P\in\Orbit_{f^{-1}}(Q)$, then~$P=(f^{-1})^n(Q)$ for some~$n\in\NN$, and hence~$f^n(P)=Q=f^{0}(Q)$. Therefore~$P\in\Orbit_f(Q)^\grand$.
\par\noindent(f)\enspace
We compute
\begin{align*}
\text{$P$ is}&\text{~not $f$-wandering} \\
&\quad\Longleftrightarrow\quad
P\in\PrePer(f) \\
&\quad\Longleftrightarrow\quad
f^n(P)=f^m(P)\quad\text{for some $n>m\ge0$,} \\
&\quad\Longleftrightarrow\quad
(f^{-1})^{m}(P)  = (f^{-1})^{n}(P)  
\quad\text{applying $(f^{-1})^{m+n}$,} \\
&\quad\Longleftrightarrow\quad
P\in\PrePer(f^{-1}) \\
&\quad\Longleftrightarrow\quad
\text{$P$ is not $f^{-1}$-wandering.}
\end{align*}
\par\noindent(g)\enspace
Two points that are grand $f$-orbit equivalent have forward orbits that differ by only finitely many points, since their forward orbits eventually merge. Hence
\[
\Bigl( \Orbit_f(P)\cup\Orbit_f(Q) \Bigr)
\setminus \Bigl( \Orbit_f(P)\cap\Orbit_f(Q) \Bigr) \quad\text{is a finite set.}
\]
Therefore~$\Orbit_f(P)$ is Zariski dense in~$X$ if and only if~$\Orbit_f(Q)$ is Zariski dense in~$X$.
\end{proof}

\subsection{Implications Between Propagation Statements}
\label{appendix:implicationspropstatements}

We prove the implications that hold universally for
the various propagation statements in Table~\textup{\ref{table:propogationstatements}}.

\begin{proposition}
\label{proposition:propagationimplications}
The following implications hold for the various orbit propagation statements in Table~\textup{\ref{table:propogationstatements}}, with the convention described in Remark~\textup{\ref{remark:(A)implies}}. These implications are illustrated in Table~\textup{\ref{table:propagationimplications}}.%
\par\noindent
\begin{minipage}{.45\hsize}
\begin{align}
\NotImply{C1A}{C1}{A} \\
\Imply{B1inftyB1}{B1\text{$\infty$}}{B1} \\
\Imply{C1inftyC1}{C1\text{$\infty$}}{C1} \\
\Imply{C2forallC1infty}{C2\text{$\forall$}}{C1\text{$\infty$}} \\
\Imply{C2forallC2exists}{C2\text{$\forall$}}{C2\text{$\exists$}} \\
\Imply{C1B1}{C1}{B1} \\
\Imply{C2existsC1}{C2\text{$\exists$}}{C1} 
\end{align}
\end{minipage}
\qquad
\vrule height 4.2\baselineskip depth 4.2\baselineskip width .5pt\relax
\begin{minipage}{.45\hsize}
\begin{align}
\Imply{C1C2exists}{C1}{C2\text{$\exists$}} \\
\Imply{C1inftyB1infty}{C1\text{$\infty$}}{B1\text{$\infty$}} \\
\Imply{C3existsA}{C3\text{$\exists$}}{A} \\
\Imply{C3forallC3exists}{C3\text{$\forall$}}{C3\text{$\exists$}} \\
\Imply{C3forallC2forall}{C3\text{$\forall$}}{C2\text{$\forall$}} \\
\Imply{C3existsC2exists}{C3\text{$\exists$}}{C2\text{$\exists$}} \\
\Imply{C1inftyC2forall}{C1\text{$\infty$}}{C2\text{$\forall$}} 
\end{align}
\end{minipage}
\end{proposition}
\begin{proof}
\NotImplyP{C1A}{C1}{A}
If the map~$f$ admits a non-trivial fibration, then~(C1) is true, since any finite set of grand $f$-orbits will be contained in a finite number of fibers, but~(A) is not true, since every grand orbits is contained in a fiber, so there are no Zariski dense grand orbits. We note that it is not clear if~(B1$\infty$) or~(C1$\infty$) implies (A), but neither do we see an obvious counterexample.
\ImplyP{B1inftyB1}{B1$\infty$}{B1}
Let~$\G_1,\ldots,\G_r$ be a collection of $f$-orbits, say $\G_i=\Orbit_f(P_i)$. Then $Y=\{P_1,\ldots,P_r\}$ is a proper Zariski closed subset of~$X$, so applying~(B1$\infty$) to~$Y$ gives~(B1).
\ImplyP{C1inftyC1}{C1$\infty$}{C1}
Let~$\G_1,\ldots,\G_r$ be a collection of grand $f$-orbits, say $\G_i=\Orbit_f^\grand(P_i)$. Then $Y=\{P_1,\ldots,P_r\}$ is a proper Zariski closed subset of~$X$, so applying~(C1$\infty$) to~$Y$ gives~(C1).
\ImplyP{C2forallC1infty}{C2$\forall$}{C1$\infty$}
Let $Y\subsetneq{X}$ be a proper Zariski closed set, and let~$\Qcal\subset{X(K)}$ be any complete set of representatives for~$X(K)/\fgeq$. Thus every grand $f$-orbit~$\G$ that contains a point of~$X(K)$ has the property that $\#(\G\cap\Qcal)=1$. We create a modified version of~$\Qcal$, which we denote by~$\Qcal'$, as follows: For each~$Q\in\Qcal$, if
\[
\Orbit_f^\grand(Q)\cap{Y(K)}\ne\emptyset,
\]
then we replace~$Q$ with one of the points in~\text{$\Orbit_f^\grand(Q)\cap{Y(K)}$}. This modified set~$\Qcal'$ is still a complete set of representatives for~$X(K)/\fgeq$. Our assumption that~(C2$\forall$) is true tells us that~$\Qcal'$ is Zariski dense in~$X$. Since~$Y$ is not Zariski dense in~$X$, it follows that
\begin{equation}
\label{eqn:QcalprimeminusY1}
\text{$\Qcal'\setminus{Y}$ is Zariski dense in~$X$.}
\end{equation}
On the other hand, the construction of~$\Qcal'$ ensures that
\begin{equation}
\label{eqn:QcalprimeminusY2}
Q\in\Qcal'\cap\bigcup_{P\in Y(K)} \Orbit_f^\grand(P)
\quad\Longrightarrow\quad
Q\in Y. 
\end{equation}
Combining~\eqref{eqn:QcalprimeminusY1} and~\eqref{eqn:QcalprimeminusY2}, we deduce that
\[
\Qcal'\setminus \biggl(\bigcup_{P\in Y(K)} \Orbit_f^\grand(P)\biggr)
\supset Q\setminus{Y}\quad\text{is Zariski dense in $X$,}
\]
and since~$\Qcal'$ is a subset of~$X(K)$, we conclude that~(C1$\infty$) is true.
\par
We remark that if we use the same argument to try to prove that (C2$\exists$) implies (C1$\infty$), we run into the problem that the replacement procedure used to create~$\Qcal'$ could potentially replace every point of~$\Qcal$ with a point in~$Y$, which prevents us from concluding that~$\Qcal'$ is Zariski dense in $X$.
\ImplyP{C2forallC2exists}{C2$\forall$}{C2$\exists$}
There certainly exists at least one set of representatives~$\Qcal\subset{X(K)}$ for
\[
X(K)/\fgeq.
\]
Our assumption that~(C2$\forall$) holds tells us that~$\Qcal$ is Zariski dense. Hence there exists at least one Zariski dense set of representatives, so~(C2$\exists$) is true.
\ImplyP{C1B1}{C1}{B1}
Let $\G_1,\ldots,\G_r$ be $f$-orbits, and let~$\L_1,\ldots,\L_s$ be the grand $f$-orbits generated by the points in~$\G_1,\ldots,\G_r$. Then each~$\G_i$ is contained in a (unique)~$\L_j$, so we have the following inclusion:
\[
\overbrace{X(K)\setminus(\L_1\cup\cdots\cup\L_s)}^{\hidewidth\text{(C1) tells us that this set is Zariski dense.}\hidewidth}
\subseteq
\underbrace{X(K)\setminus(\G_1\cup\cdots\cup\G_r) }_{\hidewidth\text{Hence this set is also Zariski dense.}\hidewidth}
\]
\ImplyP{C2existsC1}{C2\text{$\exists$}}{C1} 
We start by extending~$K$ so that~$X(K)\ne\emptyset$. The assumption that~(C2$\exists$) is true says that, after replacing~$K$ by a further finite extension, there is a set of points $\Qcal\subseteq{X(K)}$ such that~$\Qcal$ is Zariski dense in~$X$ and such that~$\Qcal$ is a complete set of representatives for~$X_f(K)/\fgeq$. 
\par
Let~$\G_1,\ldots,\G_r$ be grand $f$-orbits. Then for each~$i$, there is a unique point~\text{$Q_i\in\Qcal\cap\G_i$}. It follows that
\[
\Qcal\setminus\{Q_1,\ldots,Q_r\} \subseteq
X(K)\setminus(\G_1\cup\cdots\cup\G_r).
\]
By assumption, the set~$\Qcal$ is Zariski dense in~$X$, and thus the same is true of \text{$\Qcal\setminus\{Q_1,\ldots,Q_r\}$}, since $\dim(X)\ge1$. Hence \text{$X(K)\setminus(\G_1\cup\cdots\cup\G_r)$} contains a Zariski dense set, which completes the proof that~(C1) is true.
\ImplyP{C1C2exists}{C1}{C2\text{$\exists$}} 
We start by enumerating the countably many\footnote{The fact that a projective variety defined over a countable field has countably many Zariski closed subsets may be reduced to the same assertion for~$\PP^N$, and then it follows from the fact that there are only countably many ideals in the Noetherian ring~$K[X_0,\ldots,X_n]$.} proper Zariski closed subsets of~$X$, say they are~$Z_1,Z_2,Z_3,\ldots\,$. We select a sequence
of points~$P_1,P_2,\ldots$ in~$X(K)$ using the following algorithm:
\begin{parts}
\Part{\textbullet}
The assumption that~(C1) is true tells us in particular that~$X(K)$ is Zariski dense in~$X$, so in any case we know that~$X(K)\ne\emptyset$. So we can choose a point~$P_1\in{X(K)}$. Let $\G_1=\Orbit_f^\grand(P_1)$. 
\Part{\textbullet}
\textsf{LOOP} $i=1,2,3,4,\ldots$
\Part{\textbullet}
The assumption that~(C1) is true tells us that $X(K)\setminus(\G_1\cup\cdots\cup\G_i)$ is Zariski dense in~$X$. The union of the proper Zariski closed subsets $Z_1,\ldots,Z_i$ is not Zariski dense in~$X$, so we can choose a point~$P_{i+1}$ satisfying
\[
P_{i+1} \in X(K) \setminus \Bigl( (\G_1\cup\cdots\cup\G_i) 
\cup (Z_1\cup\cdots\cup Z_i) \Bigr).
\]
We let $\G_{i+1}=\Orbit_f^\grand(P_{i+1})$. 
\Part{\textbullet}
\textsf{END LOOP}
\end{parts}
The output from this algorithm is a set of points
\[
\Pcal = \{P_1,P_2,P_3,\ldots\}\subseteq X(K)
\]
whose grand orbits are distinct. We claim that~$\Pcal$ is Zariski dense. If not, then~$\overline{\Pcal}$ is a proper Zariski closed subset of~$X$,
so by our enumeration~$Z_1,Z_2,\ldots$ of all of the proper Zariski closed subsets of~$X$, there is an index~$j$ such that~$\Pcal\subseteq{Z_j}$. But $P_{j+1}\in\Pcal$ and, by construction, we have $P_{j+1}\notin{Z_j}$. This contradiction shows that~$\Pcal$ is a Zariski dense set of points in~$X(K)$ that represent distinct elements of~$X(K)/\fgeq$. Extending~$\Pcal$ to a full set of representatives shows that~(C2$\exists$) is true.
\ImplyP{C1inftyB1infty}{C1$\infty$}{B1$\infty$}
We prove these simultaneously.
Every $f$-orbit is contained in a (unique) grand $f$-orbit, since for every point~$P$ we have
\[
\Orbit_f(P) \subseteq \Orbit_f^\grand(P).
\]
This implies that if~$\Pcal\subset{X(K)}$ is any set of points, then
\[
\biggl( X(K)\setminus\bigcup_{P\in\Pcal} \Orbit_f(P) \biggr)
\supseteq \biggl( X(K)\setminus\bigcup_{P\in\Pcal} \Orbit_f^\grand(P) \biggr).
\]
Taking~$\Pcal$ to be a finite set of points gives \text{(C1)}$\;\longrightarrow\;$\text{(B1)}, and taking~$\Pcal=Y(K)$ gives \text{(C1$\infty$)}$\;\longrightarrow\;$\text{(B1$\infty$)}.
\ImplyP{C3existsA}{C3$\exists$}{A}
We start by extending~$K$ so that~$X(K)\ne\emptyset$. The assumption that~(C3$\exists$) is true says that, after replacing~$K$ by a further finite extension, there is a set of points $\Qcal\subseteq{X(K)}$ such that~$\Qcal$ is Zariski dense in~$X$ and such that~$\Qcal$ is a complete set of representatives for~$X_f^\dense(K)/\fgeq$. Let~$Q\in\Qcal$. Lemma~\ref{lemma:orbitelemproperties}(j) tells us that~$Q\in{X_f^\dense}$, and we know~$Q\in{X(K)}$, so we have found a point~$Q\in{X_f^\dense(K)}$. By definition, the orbit~$\Orbit_f(Q)$ is Zariski dense in~$X$, which proves that~(A) is true.
\ImplyP{C3forallC3exists}{C3$\forall$}{C3$\exists$}
By assumption, there is a point $Q\in{X_f^\dense}(K)$. It follows that
$X_f^\dense(K)/\fgeq$ is non-empty, so there exists a set~$\Qcal\subset{X_f^\dense(K)}$ of  representatives for $X_f^\dense(K)/\fgeq$. The assumption~(C3$\forall$) tells us that~$\Qcal$ is Zariski dense in~$X$, so the existence of~$\Qcal$ implies that~(C3$\exists$) is true.
\ImplyP{C3forallC2forall}{C3$\forall$}{C2$\forall$}
Let~$\Qcal\subset{X(K)}$ be a complete set of representatives for $X(K)/\fgeq$. Then $\Qcal\cap{X^\dense_f(K)}$ is a complete set of representative for~$X_f^\dense(K)/\fgeq$, where we are using Lemma~\ref{lemma:orbitelemproperties}(j). The assumption that~(v3$\forall$) is true tells us that~$\Qcal\cap{X^\dense_f(K)}$ is Zariski dense in~$X$, from which it is clear that~$\Qcal$ is Zariski dense in~$X$.
\ImplyP{C3existsC2exists}{C3$\exists$}{C2$\exists$}
We are given that there exists a set~$\Qcal\subset{X_f^\dense(K)}$ such that~$\Qcal$ is Zariski dense in~$X$ and such that~$\Qcal$ is a complete set of representatives for~$X_f^\dense(K)/\fgeq$. Let~$\Qcal'\subset{(X\setminus{X}_f^\dense)(K)}$ be a complete set of representatives for
\[
(X\setminus X_f^\dense)(K)/\fgeq.
\]
Lemma~\ref{lemma:orbitelemproperties}(j) tells us that~$\Qcal\cup\Qcal'$ is a complete set of representatives for~$X(K)/\fgeq$, and that~$\Qcal\cup\Qcal'$ is Zariski dense in~$X$, since it contains~$\Qcal$. Hence~(C2$\exists$) is true.
\ImplyP{C1inftyC2forall}{C1$\infty$}{C2$\forall$} 
Let~$\Qcal\subset{X(K)}$ be a complete set of representatives for~$X(K)/\fgeq$, and let
\[
Y = \overline{\Qcal} = \text{the Zariski closure of $\Qcal$.}
\]
We suppose that~$Y\ne{X}$, and our goal is to show that~(C1$\infty$) is false for this~$Y$. The choice of~$\Qcal$ tells us that
\begin{equation}
\label{eqn:XKcupQQcalgrand}
X(K) = \bigcup_{Q\in\Qcal} \Bigl( \Orbit_f^\grand(Q)\cap X(K) \Bigr).
\end{equation}
(Indeed, it even tells us that the union is a disjoint union, but we will not need this fact.) The facts~$Y=\overline\Qcal$ and~$\Qcal\subseteq{X(K)}$ imply that
\begin{equation}
\label{eqn:XksetminuscupYQcal}
Y(K)\supseteq\Qcal.
\end{equation}
Hence
\begin{align}
\label{eqn:XknegPinYKgrandorb}
X(K) \setminus \biggl( \bigcup_{P\in Y(K)} & \Orbit_f^\grand(P) \biggr) \\
&\subseteq X(K) \setminus 
\smash[t]{ \biggl( \bigcup_{Q\in\Qcal} \Orbit_f^\grand(Q) \biggr) }
\quad\text{from \eqref{eqn:XksetminuscupYQcal},} \notag \\
&= \emptyset \quad\text{from \eqref{eqn:XKcupQQcalgrand}.} \notag
\end{align}
This certainly contradicts~(C1$\infty$), since~(C1$\infty$) would imply that the set~\eqref{eqn:XknegPinYKgrandorb} is Zariski dense in~$X$.
\end{proof}

The next result describes another relation among the various propagation properties, but it has a slightly different flavor from the results in Proposition~\ref{proposition:propagationimplications}. It also probably does not have wide applicability, since for example it cannot be used if the $f$-preperiodic points are Zariski dense; but it is exactly what we need for linear maps of~$\PP^N$.

\begin{lemma}
\label{lemma:C2forallZDimpliesC3forall}
Let $f:X\to{X}$.  Then
\[
\text{\textup{(A)}\quad and\quad \textup{(C2$\forall$)}\quad
and\quad 
$\left[
\text{\TAB{$X_f^\dense$ contains
a non-\\
empty Zariski open set\\}}
\right]$
}
\quad\Longrightarrow\quad
\textup{(C3$\forall)$}.
\]
\end{lemma}
\begin{proof}
Let
\[
\Qcal\subset X_f^\dense(K)
\quad\text{such that}\quad
\Qcal \xleftrightarrow{\;\text{bijective}\;}
X_f^\dense(K)/\fgeq.
\]
Our goal is to prove that~$\Qcal$ is Zariski dense in~$X$. To do this, we consider the complement of~$X_f^\dense$, and we choose a set
\[
\Qcal'\subset (X\setminus X_f^\dense)(K)
\quad\text{such that}\quad
\Qcal' \xleftrightarrow{\;\text{bijective}\;}
(X\setminus X_f^\dense)(K)/\fgeq.
\]
Then the union satisfies
\[
\Qcal\cup\Qcal'
\xleftrightarrow{\;\text{bijective}\;} X(K)/\fgeq,
\]
so the assumption that~(C2$\forall$) is true tells us that~$\Qcal\cup\Qcal'$ is Zariski dense in~$X$. However, the set~$\Qcal'$ is not Zariski dense in~$X$, since
\[
\Qcal' \subset X\setminus X_f^\dense,
\]
and the assumption that~$X_f^\dense$ contains a non-empty Zariski open set implies that the complement of~$X_f^\dense$ is not Zariski dense. Hence~$\Qcal$ is Zariski dense in~$X$, which completes the proof that~$f$ has propagation property~(C3$\forall$).
\end{proof}

We will apply the next lemma, which is undoubtedly well known, to the case of an automorphism. But since it is no harder to handle the case of arbitrary backward branches, we formulate it in that generality.

\begin{lemma}
\label{lemma:zariskidensebranch}
Let~$P_0\in{X}$ be a point whose forward orbit~$\Orbit_f(P_0)$ is Zariski dense in~$X$, and let
\[
\Bcal = \{P_0,P_1,P_2,\ldots\}\subset X(\Kbar)
\]
be a complete backward $f$-branch of~$P_0$, i.e., a sequence of points in~$X(\Kbar)$ satisfying
\[
f(P_{i+1}) = P_i\quad\text{for all $i\ge0$.}
\]
Then~$\Bcal$ is Zariski dense in~$X$.
\end{lemma}
\begin{proof}
Our first observation is that~$P_0,P_1,\ldots$ are distinct. To see why, suppose that~$P_i=P_j$ for some $i>j$. Then
\[
f^{i-j}(P_0) = f^{i-j}\bigl(f^j(P_j)\bigr) 
= \underbrace{f^i(P_j) = f^i(P_i)}_{\text{since $P_i=P_j$}} = P_0,
\]
which implies that~$P_0$ is periodic and thus contradicts the assumed Zariski density of~$\Orbit_f(P_0)$. (Note we always assume that $\dim(X)\ge1$.)
\par
Our goal is to prove that the Zariski closure~$\overline\Bcal$ of~$\Bcal$ is equal to~$X$. 
For each integer~$m\ge0$, we let
\[
\Bcal_m := \{P_m,P_{m+1},P_{m+2},\ldots\}
\]
denote the branch with it's first~$m$ points omitted. The actions of~$f$ on these sets and their closures satisfy
\begin{equation}
\label{eqn:fBmplus1toBm}
f : \Bcal_{m+1} \longrightarrow \Bcal_m
\quad\text{and}\quad
f : \overline{\Bcal_{m+1}} \longrightarrow \overline{\Bcal_m}.
\end{equation}
\par
We now suppose that
\begin{equation}
\label{eqn:PmnotinBmplus1}
P_m\notin\overline{\Bcal_{m+1}}\quad\text{for all $m\ge0$,}
\end{equation}
and derive a contradiction. Since~$\Bcal_m=\{P_m\}\cup\Bcal_{m+1}$, it would follow from~\eqref{eqn:PmnotinBmplus1} that the single point~$\{P_m\}$ is an irreducible component of~$\overline{\Bcal_m}$; and hence, since~$\Bcal\setminus\Bcal_m$ is a finite set, it would follow that~$\{P_m\}$ is an irreducible component of~$\overline\Bcal$. The Zariski closed set~$\overline\Bcal$ has only finitely many irreducible components, and we proved earlier that the~$P_i$ are distinct, which completes our proof that~\eqref{eqn:PmnotinBmplus1} is false. 
\par
Therefore, there exists an index~$\ell\ge0$ such that
\begin{equation}
\label{eqn:PellinZCBellplus1}
P_\ell \in \overline{\Bcal_{\ell+1}}.
\end{equation}
This implies that
\[
\overline{\Bcal_{\ell}} 
= \overline{\{P_\ell\}\cup\Bcal_{\ell+1}}
= \{P_\ell\}\cup\overline{\Bcal_{\ell+1}}
= \overline{\Bcal_{\ell+1}},
\]
and then~\eqref{eqn:fBmplus1toBm} tells us that~$f$ induces a map
\begin{equation}
\label{eqn:fBellplus1toitself}
f : \overline{\Bcal_{\ell+1}} \longrightarrow \overline{\Bcal_{\ell+1}}.
\end{equation}
We know from~\eqref{eqn:PellinZCBellplus1} that the point~$P_\ell$ is in~$\Bcal_{\ell+1}$, so~\eqref{eqn:fBellplus1toitself} implies that
\begin{equation}
\label{eqn:orbitfPell}
\Orbit_f(P_\ell) \subset \overline{\Bcal_{\ell+1}}.
\end{equation}
We compute
\begin{align*}
X &= \overline{\Orbit_f(P_0)} 
&&\text{by assumption,} \\
&\subseteq \overline{\Orbit_f(P_\ell)} 
&&\text{since $f^\ell(P_\ell)=P_0$,} \\
&\subseteq \overline{\Bcal_{\ell+1}}
&&\text{from \eqref{eqn:orbitfPell},} \\
&= \overline{ \Bcal \setminus \{P_0,\ldots,P_{\ell}\} } 
&&\text{by definition of $\Bcal_m$,} \\
&\subseteq \overline\Bcal.
\end{align*}
Hence~$\overline\Bcal=X$, which completes the proof that~$\Bcal$ is Zariski dense in~$X$.
\end{proof}

If~$f:X\to{X}$ is an automorphism, then the existence of a Zariski dense orbit automatically implies a weak propagation property, as in the following result. The proof follows from the fact that~$\Orbit_f(P)$ is Zariski dense if and only if~$\Orbit_{f^{-1}}(P)$ is Zariski dense. This allows us to use points in the backward orbit to prove that~(A) implies propagation property~(B1); but this is a bit of a cheat, and in any case the idea cannot be used to prove a grand orbit propagation property such as~(C1).

\begin{proposition}
\label{proposition:AimpliesB1forautomorphisms}
Let $f:X\to{X}$ be an automorphism. Then
\[
\textup{(A)} \quad\Longrightarrow\quad \textup{(B1)}
\]
\end{proposition}
\begin{proof}
We are given that there is a point~$P_0\in{X(K)}$ whose $f$-orbit $\Orbit_f(P_0)$ is Zariski dense. The forward~$f^{-1}$-orbit of~$P_0$ is a backward $f$-branch of~$P_0$ in the sense of Lemma~\ref{lemma:zariskidensebranch}, so Lemma~\ref{lemma:zariskidensebranch} tells us that~$\Orbit_{f^{-1}(P_0)}$ is Zariski dense in~$X$.
\par
We now commence the proof of propagation property~(B1). Let $\G_1,\ldots,\G_r$ be a collection of~$f$-orbits. The fact that the~$\G_i$ are forward $f$-orbits implies that they have only finitely many points in common with the $f^{-1}$-orbit~$\Orbit_{f^{-1}}(P_0)$.\footnote{Suppose that~$\Orbit_f(Q)\cap\Orbit_{f^{-1}}(P_0)\ne\emptyset$. Then~$Q=f^{-m}(P_0)$ for some~$m\ge0$, so $\Orbit_f(Q)\cap\Orbit_{f^{-1}}(P_0)=\bigl\{f^n(P_0):-m\le{n}\le0\bigr\}$ is a finite set.} Hence
\begin{align*}
X(K)\setminus(\G_1\cup\cdots\cup\G_r) 
&\supseteq \Orbit_{f^{-1}}(P_0) \setminus (\G_1\cup\cdots\cup\G_r) \\
&= \underbrace{\Orbit_{f^{-1}}(P_0) \setminus \{\text{finite set}\}}_{\text{Zariski dense in $X$.}}.
\end{align*}
This completes the proof that~$f$ has propagation property~(B1).
\end{proof}

\begin{remark}
We note that if we assume the dynamical Mordell--Lang conjecture for~$f^{-1}$, then we can strengthen Proposition~\ref{proposition:AimpliesB1forautomorphisms} to the statement that~(A) implies that~$f$ has propagation property~(B1$\infty$). Briefly, if a Zariski closed subset~$Y\subseteq{X}$ contains~$f^{-n}(P_0)$ for infinitely many~$n\ge0$, then Mordell--Lang says that~$Y$ contains an~$f^{-1}$-invariant subvariety~$Z$ that contains those points. In particular, there is a point~$f^{-m}(P_0)\in{Z}$. But then~$Z$ is also~$f$-invariant, so it contains~$f^m(f^{-m}(P_0))=P_0$. Hence~$\Orbit_f(P_0)\subset{Z}$, and then the Zariski density of~$\Orbit_f(P_0)$ forces~$Z=X$, and thus also $Y=X$.
\end{remark}

\subsection{A Coarse Height Counting Estimate}
\label{appendix:htcountPN}

There are many stronger results that imply the following height counting lemma, but we include it to indicate that the estimate that we need can be obtained in an elementary manner.

\begin{lemma}
\label{lemma:ctYKTleCTN}
Let~$K$ be a number field, and let~$N\ge1$, and compute counting functions using the multiplicative $K$-height on $\PP^N(K)$; see~\textup{\cite[Section~B.2]{MR1745599}}.
\begin{parts}
\Part{(a)}
There are constants~$\Cl{PNapp1}(K,N)>0$ and~$\Cl{PNapp2}(K,N)>0$ such that
\[
\Cr{PNapp1}(K,N) T^{N+1} \le \Count\bigl(\PP^N(K),T\bigr) \le \Cr{PNapp2}(K,N) T^{N+1}.
\]
\Part{(b)}
Let $Y\subsetneq\PP^N$ be a Zariski closed set. 
There is a constant~$\Cl{YinPNlemma}(K,Y)$ such that
\begin{equation}
\label{eqn:CtYKTCTdimy1}
\Count\bigl(Y(K),T\bigr) \le \Cr{YinPNlemma}(K,Y)T^{1+\dim(Y)},
\end{equation}
where the dimension of an algebraic subset of~$\PP^N$ is defined to be the maximum of the dimensions of its geometrically irreducible components. In particular, if~$Y\subsetneq\PP^N$ is a proper subset, then
\[
\Count\bigl(Y(K),T\bigr) \le \Cr{YinPNlemma}(K,Y)T^N.
\]
\end{parts}
\end{lemma}
\begin{proof}
(a)\enspace
Schanuel's formula~\cite{MR557080} gives a precise asymptotic formula for~$\Count\bigl(\PP^N(K),T\bigr)$, but the weak estimate that we have cited in~(a) follows from the standard proofs that $\PP^N(K)$ has only finitely many points of bounded height.
\par\noindent(b)\enspace
We prove~\eqref{eqn:CtYKTCTdimy1} by induction on~$\dim(Y)$.
If~$\dim(Y)=0$, then~$Y$ is a finite set of points, so its counting function is bounded and~\eqref{eqn:CtYKTCTdimy1} is trivially true.
\par
Suppose now that~$m\ge1$, that we know~\eqref{eqn:CtYKTCTdimy1} for all algebraic sets of dimension at most~$m-1$, and that~$Y\subseteq\PP^N$ is a Zariski closed subset with~$\dim(Y)=m$. (Of course, we must have $m\le{N}$.) Let~$Y_1,\ldots,Y_r$ be the irreducible components of~$Y$ of dimension~$m$. By the induction hypothesis, it suffices to prove~\eqref{eqn:CtYKTCTdimy1} for the union of the~$Y_i$, which reduces us to the case that~$Y$ is irreducible and of dimension~$m$. A generic projection of~$Y$ onto a linear subspace of dimension~$m$ yields a quasi-finite rational map
\[
\f : Y \dashrightarrow \PP^m.
\]
Let~$Y^\circ\subseteq{Y}$ be a non-empty Zariski open subset on which~$\f$ is a quasi-finite morphism of degree~$d$. Then
\begin{align}
\label{eqn:CtYoKT}
\Count\bigl( Y^\circ(K), T\bigr) 
&\le (\deg\f)\cdot\Count\bigl(\PP^m(K),T\bigr)\notag\\
&\le (\deg\f)\cdot \Cl{CtYoKT}(m,K) T^{1+m} quad\text{from (a)}\notag\\
&= \Cr{CtYoKT}(K,Y) T^{1+\dim(Y)}.
\end{align}
We next observe that~$Y\setminus{Y^\circ}$ is a Zariski closed subset of~$\PP^N$ of dimension at most~$m-1$, so our induction hypothesis gives
\begin{equation}
\label{eqn:CtYminusYoKT}
\Count\bigl( (Y\setminus Y^\circ)(K),T\bigr) 
\le \Cl{CtYminusYoKT}(K,Y) T^{1+\dim(Y\setminus Y^\circ)}
\le \Cr{CtYminusYoKT}(K,Y) T^{\dim(Y)}.
\end{equation}
Combining~\eqref{eqn:CtYoKT} and~\eqref{eqn:CtYminusYoKT} yields~\eqref{eqn:CtYKTCTdimy1} for~$Y$, which completes our induction proof that~\eqref{eqn:CtYKTCTdimy1} is true for all algebraic subsets of~$\PP^N$. In particular, if~$Y\subsetneq\PP^N$ is a proper Zariski closed subset, then its dimension is at most~$N-1$. This completes the proof of Lemma~\ref{lemma:ctYKTleCTN}(b).
\end{proof}




\begin{thebibliography}{10}

\bibitem{MR2862064}
E.~Amerik, F.~Bogomolov, and M.~Rovinsky.
\newblock Remarks on endomorphisms and rational points.
\newblock {\em Compos. Math.}, 147(6):1819--1842, 2011.

\bibitem{MR2784670}
Ekaterina Amerik.
\newblock Existence of non-preperiodic algebraic points for a rational self-map
  of infinite order.
\newblock {\em Math. Res. Lett.}, 18(2):251--256, 2011.

\bibitem{MR2766180}
J.~P. Bell, D.~Ghioca, and T.~J. Tucker.
\newblock The dynamical {M}ordell-{L}ang problem for \'{e}tale maps.
\newblock {\em Amer. J. Math.}, 132(6):1655--1675, 2010.
\newblock MR2766180.

\bibitem{MR1255693}
Gregory~S. Call and Joseph~H. Silverman.
\newblock Canonical heights on varieties with morphisms.
\newblock {\em Compositio Math.}, 89(2):163--205, 1993.
\newblock MR1255693.

\bibitem{MR2597304}
Thomas Dedieu.
\newblock Severi varieties and self-rational maps of {$K3$} surfaces.
\newblock {\em Internat. J. Math.}, 20(12):1455--1477, 2009.

\bibitem{MR1109353}
Gerd Faltings.
\newblock Diophantine approximation on abelian varieties.
\newblock {\em Ann. of Math. (2)}, 133(3):549--576, 1991.

\bibitem{MR337997}
Gerhard Frey and Moshe Jarden.
\newblock Approximation theory and the rank of abelian varieties over large
  algebraic fields.
\newblock {\em Proc. London Math. Soc. (3)}, 28:112--128, 1974.

\bibitem{MR1183561}
Tsuyoshi Fujiwara.
\newblock Varieties of small {K}odaira dimension whose cotangent bundles are
  semiample.
\newblock {\em Compositio Math.}, 84(1):43--52, 1992.

\bibitem{MR3829742}
Dragos Ghioca and Fei Hu.
\newblock Density of orbits of endomorphisms of commutative linear algebraic
  groups.
\newblock {\em New York J. Math.}, 24:375--388, 2018.

\bibitem{MR4396628}
Dragos Ghioca and Sina Saleh.
\newblock Zariski dense orbits for regular self-maps on split semiabelian
  varieties.
\newblock {\em Canad. Math. Bull.}, 65(1):116--122, 2022.

\bibitem{MR3937327}
Dragos Ghioca and Matthew Satriano.
\newblock Density of orbits of dominant regular self-maps of semiabelian
  varieties.
\newblock {\em Trans. Amer. Math. Soc.}, 371(9):6341--6358, 2019.

\bibitem{MR3557780}
Dragos Ghioca and Thomas Scanlon.
\newblock Density of orbits of endomorphisms of abelian varieties.
\newblock {\em Trans. Amer. Math. Soc.}, 369(1):447--466, 2017.

\bibitem{MR0463157}
Robin Hartshorne.
\newblock {\em Algebraic geometry}.
\newblock Graduate Texts in Mathematics, No. 52. Springer-Verlag, New
  York-Heidelberg, 1977.

\bibitem{MR1745599}
Marc Hindry and Joseph~H. Silverman.
\newblock {\em Diophantine geometry}, volume 201 of {\em Graduate Texts in
  Mathematics}.
\newblock Springer-Verlag, New York, 2000.
\newblock An introduction.

\bibitem{MR3586372}
Daniel Huybrechts.
\newblock {\em Lectures on {K}3 surfaces}, volume 158 of {\em Cambridge Studies
  in Advanced Mathematics}.
\newblock Cambridge University Press, Cambridge, 2016.
\newblock MR3586372.

\bibitem{MR4530188}
Jia Jia, Junyi Xie, and De-Qi Zhang.
\newblock Surjective endomorphisms of projective surfaces: the existence of
  infinitely many dense orbits.
\newblock {\em Math. Z.}, 303(2):Paper No. 39, 23, 2023.

\bibitem{MR3456169}
Shu Kawaguchi and Joseph~H. Silverman.
\newblock Dynamical canonical heights for {J}ordan blocks, arithmetic degrees
  of orbits, and nef canonical heights on abelian varieties.
\newblock {\em Trans. Amer. Math. Soc.}, 368(7):5009--5035, 2016.
\newblock MR3456169.

\bibitem{MR65218}
Serge Lang and Andr\'{e} Weil.
\newblock Number of points of varieties in finite fields.
\newblock {\em Amer. J. Math.}, 76:819--827, 1954.

\bibitem{MR2917181}
Jun Li and Christian Liedtke.
\newblock Rational curves on {K}3 surfaces.
\newblock {\em Invent. Math.}, 188(3):713--727, 2012.
\newblock MR2917181.

\bibitem{MR3871505}
Yohsuke Matsuzawa, Kaoru Sano, and Takahiro Shibata.
\newblock Arithmetic degrees and dynamical degrees of endomorphisms on
  surfaces.
\newblock {\em Algebra Number Theory}, 12(7):1635--1657, 2018.
\newblock MR3871505.

\bibitem{MR1323985}
Michael McQuillan.
\newblock Division points on semi-abelian varieties.
\newblock {\em Invent. Math.}, 120(1):143--159, 1995.

\bibitem{MR3126567}
Alice Medvedev and Thomas Scanlon.
\newblock Invariant varieties for polynomial dynamical systems.
\newblock {\em Ann. of Math. (2)}, 179(1):81--177, 2014.

\bibitem{mengzhang2023}
Sheng Meng and De-Qi Zhang.
\newblock Advances in the equivariant minimal model program and their
  applications in complex and arithmetic dynamics, 2023.
\newblock \url{arXiv:2311.16369}.

\bibitem{MR2514037}
David Mumford.
\newblock {\em Abelian varieties}, volume~5 of {\em Tata Institute of
  Fundamental Research Studies in Mathematics}.
\newblock Published for the Tata Institute of Fundamental Research, Bombay; by
  Hindustan Book Agency, New Delhi, 2008.
\newblock With appendices by C. P. Ramanujam and Yuri Manin, Corrected reprint
  of the second (1974) edition, MR2514037.

\bibitem{MR2551469}
Noboru Nakayama and De-Qi Zhang.
\newblock Building blocks of \'{e}tale endomorphisms of complex projective
  manifolds.
\newblock {\em Proc. Lond. Math. Soc. (3)}, 99(3):725--756, 2009.

\bibitem{MathOverflow267708}
Francesco Polizzi.
\newblock Infinitely many exceptional curves on ruled surfaces.
\newblock MathOverflow.
\newblock \url{https://mathoverflow.net/q/267708} (version: 2020-06-15).

\bibitem{MR3210707}
Bjorn Poonen.
\newblock {$p$}-adic interpolation of iterates.
\newblock {\em Bull. Lond. Math. Soc.}, 46(3):525--527, 2014.

\bibitem{MR557080}
Stephen~Hoel Schanuel.
\newblock Heights in number fields.
\newblock {\em Bull. Soc. Math. France}, 107(4):433--449, 1979.

\bibitem{MSE493989}
Brian~M. Scott.
\newblock Difference of closures and closure of difference.
\newblock Mathematics Stack Exchange.
\newblock \url{https://math.stackexchange.com/q/493989} (version: 2013-09-15).

\bibitem{MR2179691}
Jean-Pierre Serre.
\newblock {\em Lie algebras and {L}ie groups}, volume 1500 of {\em Lecture
  Notes in Mathematics}.
\newblock Springer-Verlag, Berlin, 2006.
\newblock 1964 lectures given at Harvard University, Corrected fifth printing
  of the second (1992) edition.

\bibitem{MR1009803}
Joseph~H. Silverman.
\newblock Integral points on curves and surfaces.
\newblock In {\em Number theory (Ulm, 1987)}, volume 1380 of {\em Lecture Notes
  in Math.}, pages 202--241. Springer, New York, 1989.

\bibitem{MR2793038}
Michel Waldschmidt.
\newblock On the {$p$}-adic closure of a subgroup of rational points on an
  {A}belian variety.
\newblock {\em Afr. Mat.}, 22(1):79--89, 2011.

\bibitem{MR3705271}
Junyi Xie.
\newblock The existence of {Z}ariski dense orbits for polynomial endomorphisms
  of the affine plane.
\newblock {\em Compos. Math.}, 153(8):1658--1672, 2017.

\bibitem{Xie2022AdelicZariskiDensity}
Junyi Xie.
\newblock The existence of zariski dense orbits for endomorphisms of projective
  surfaces.
\newblock {\em Journal of the American Mathematical Society}, 2022.
\newblock with an appendix in collaboration with T. Tucker. Electronically
  published April 18, 2022. \url{https://doi.org/10.1090/jams/1004}.

\bibitem{MR2408228}
Shou-Wu Zhang.
\newblock Distributions in algebraic dynamics.
\newblock In {\em Surveys in differential geometry. {V}ol. {X}}, volume~10 of
  {\em Surv. Differ. Geom.}, pages 381--430. Int. Press, Somerville, MA, 2006.

\end{thebibliography}

\end{document}